\documentclass[11pt]{article}
\usepackage{CJK}
\usepackage{graphicx}
\usepackage{fancyhdr}
\usepackage{amsmath,amsthm,amssymb}
\usepackage{bm}
\usepackage{mathrsfs}
\usepackage{multirow}
\usepackage{indentfirst}
\usepackage{url}
\usepackage{stmaryrd}
\usepackage{color}
\usepackage{epstopdf}
\usepackage{natbib}
\usepackage[caption=false]{subfig}
\usepackage[backref]{hyperref}
\newtheorem{theorem}{Theorem}[section]
\newtheorem{prop}{Proposition}%[chapter]
\newtheorem{remark}{Remark}

\numberwithin{equation}{section}
\makeatletter
\newcommand{\Rmnum}[1]{\expandafter\@slowromancap\romannumeral #1@}
\makeatother

\allowdisplaybreaks
\def\hDash{\bot\!\!\!\bot}

\begin{document}

\title{\bf Weighted residual empirical processes, martingale transformations, and model specification tests for regressions with diverging number of parameters
%\footnote{Lixing Zhu is a Chair professor of Department of Mathematics at Hong Kong Baptist University, Hong Kong, China. He was supported by a grant from the University Grants Council of Hong Kong, Hong Kong, China. }
}
%The authors thank the editor, the associate editor and two referees for their constructive suggestions that led to the  improvement of an early manuscript.}}
\author{Falong Tan$^1$, Xu Guo$^2$ and Lixing Zhu$^{3}$ \\~\\
{\small {\small {\it $^1$ Department of Statistics and Data Science, Hunan University, Changsha, China} }}\\
{\small {\small {\it $^2$ School of Statistics, Beijing Normal University, Beijing, China} }}\\
{\small {\small {\it $^3$ Department of Statistics, Beijing Normal University at Zhuhai, Zhuhai, China} }}
}
\date{}
\maketitle

\begin{abstract}
This paper explores hypothesis testing for the parametric forms of the mean and variance functions in regression models under diverging-dimension settings. To mitigate the curse of dimensionality, we introduce weighted residual empirical process-based tests, both with and without martingale transformations. The asymptotic properties of these tests are derived from the behavior of weighted residual empirical processes and their martingale transformations under the null and alternative hypotheses.
The proposed tests without martingale transformations achieve the fastest possible rate of detecting local alternatives, specifically of order $n^{-1/2}$, which is unaffected by dimensionality. However, these tests are not asymptotically distribution-free. To address this limitation, we propose a smooth residual bootstrap approximation and establish its validity in diverging-dimension settings.
In contrast, tests incorporating martingale transformations are asymptotically distribution-free but exhibit an unexpected limitation: they can only detect local alternatives converging to the null at a much slower rate of order $n^{-1/4}$, which remains independent of dimensionality. This finding reveals a theoretical advantage in the power of tests based on weighted residual empirical process without martingale transformations over their martingale-transformed counterparts, challenging the conventional wisdom of existing asymptotically distribution-free tests based on martingale transformations.
To validate our approach, we conduct simulation studies and apply the proposed tests to a real-world dataset, demonstrating their practical effectiveness.
\\

{\bf Key words:}  Asymptotically distribution-free, diverging number of parameters, martingale transformation, weighted residual empirical processes.
\end{abstract}

%\newpage
%\baselineskip=16pt

%\newpage

%\setcounter{equation}{0}
\section{Introduction}
This research is motivated by model checking for the mean and variance functions in high dimension settings where the number of parameters diverges as the sample size tends to infinity. Consider the following regression model:
\begin{equation}\label{1.1}
Y=m(X)+ \varepsilon, \quad \varepsilon=\sigma(X) \eta,
\end{equation}
where $(Y, X)$ is a random vector with real-valued response variable $Y$ and $d$-dimensional predictor vector $X$, $m(x)=E(Y|X=x)$ is the regression function, $\varepsilon=\sigma(X) \eta$ is the error term, $\eta$ is independent of $X$ with mean $0$ and variance $1$, and $\sigma^2(X)=var(Y|X)$ is the unknown conditional variance function. Our interest is to check whether the mean function $m(X)$ belongs to some parametric class of functions $\mathcal{M}=\{ m(\cdot, \beta): \beta \in \Theta \subset \mathbb{R}^p\}$, and whether the variance function $\sigma^2(X)$ belongs to some parametric class of variance functions $\tilde{\mathcal{M}}=\{ \sigma^2(\cdot, \theta): \theta \in \tilde{\Theta} \subset \mathbb{R}^q\}$ when the dimensions $d$, $p$ and $q$ of the predictor vector $X$, parameter vectors $\beta$ and $\theta$, respectively, diverge as the sample size $n$ tends to infinity.
The dependence of quantities $X, m, d, p, q$, and $\sigma$ on $n$ is suppressed for notational simplicity throughout this paper.

Since the pioneering works of \cite{Bierens1982} and \cite{stute1997}, the cusum processes of the residuals have formed the foundation for constructing test statistics in regression model checking. These tests are well known for their ability to detect local alternatives at the fastest possible rate of convergence in hypothesis testing. However, they are often not asymptotically distribution-free, as their limiting null distributions depend on the unknown Data Generating Process (DGP). To address this, resampling methods, such as the wild bootstrap, are commonly employed to determine critical values (see, e.g., \cite{stute1998a}, \cite{dominguez2005}). Additionally, \cite{stute1997} explored the principal components of the Cram\'{e}r-von Mises test statistic, providing a framework to approximate the limiting null distributions of such tests.
In the context of univariate predictor, \cite{stute1998b} introduced the martingale transformation for residual-marked empirical processes, enabling the derivation of tests with tractable limiting null distributions. This innovation, inspired by the martingale transformation proposed by \cite{khmaladze1981} for goodness-of-fit tests of cumulative distribution functions, has become an important methodology for constructing asymptotically distribution-free tests in model checking. In econometrics, \cite{bai2003} was the first to  apply this transformation for model specification testing. Numerous follow-up studies have extended this methodology to various testing problems, including \cite{koul1999}, \cite{bai2001}, \cite{stute2002}, \cite{koenker2002, koenker2006}, \cite{delgado2008}, \cite{khmaladze2009}, \cite{tan2019a}, and \cite{lu2020}.
Another class of tests for model checking is based on local smoothing techniques, which rely on nonparametric regression estimation. Examples include the works of \cite{hardle1993}, \cite{fanli1996}, \cite{zheng1996}, \cite{horowitz2001}, and \cite{lavergne2008, lavergne2012}. In fixed-dimension settings, these tests can be asymptotically distribution-free and are particularly sensitive to high-frequency alternative models.

However, most existing tests from these two classes in the literature are susceptible to the curse of dimensionality, as data sparsity in high-dimensional spaces poses significant challenges. For instance, as highlighted by \cite{escanciano2006a}, the residual-marked empirical process proposed by \cite{stute1997} uses the indicator function of the predictor vector $X$ that is a product of indicator functions of every component of $X$, and is highly sensitive to the dimension of multivariate predictors due to data sparsity in high-dimensional spaces. Consequently, the corresponding tests suffer severely from the curse of dimensionality.
While local smoothing-based tests (e.g., \cite{hardle1993}, \cite{zheng1996}) detect local alternatives at slower rates than empirical process-based tests, typically at the rate of order $1/\sqrt {nh^{d/2}}$ when kernel estimation is applied, where $h$ is the bandwidth. These sensitivity rates can decrease dramatically with increasing dimensionality, see also the comment in \cite{guo2016}. Therefore, \cite{escanciano2006a} and \cite{guo2016} proposed dimension reduction techniques to mitigate the curse of dimensionality in fixed-dimension settings. In diverging-dimension scenarios, however, there are limited test procedures designed to address this challenge. \cite{tan2019a, tan2022} introduced methods to handle dimensionality issues for testing single-index and multi-index models, respectively. \cite{tan2022} demonstrated that tests in diverging-dimension settings may exhibit fundamentally different limiting properties compared to those in fixed-dimension settings.
Despite these advancements, existing methods often rely critically on dimension reduction structures under the null hypothesis and cannot be easily extended to test general parametric regression models without such structures in diverging-dimension settings.

In this paper, we aim to develop model specification tests for general parametric regression models in diverging-dimension scenarios. To address the dimensionality problem, we construct the test statistic based on weighted residual empirical processes, which rely on one-dimensional error terms instead of high-dimensional predictor vector. This methodology builds on the foundational work of \cite{stute2008} and \cite{escanciano2010} in fixed-dimension settings.
By focusing on the univariate error terms, the corresponding tests circumvent data sparsity issues in high-dimensional spaces, thereby significantly alleviating the curse of dimensionality. Further discussion on this issue is provided in Section~2.
Leveraging high dimensional empirical process theory, we investigate the asymptotic properties of the weighted residual empirical process under both the null and alternative hypotheses, provided that $p=o(n^{1/3}/\log^2{n})$. We demonstrate that tests based on this process can detect local alternatives distinct from the null at the parametric rate of order $n^{-1/2}$.
Since these tests are not asymptotically distribution-free, we propose a smoothed residual bootstrap to approximate their limiting null distribution and establish the asymptotic validity of the approximation under diverging-dimension settings.

To reduce the computational burden of bootstrapping, we propose a martingale transformation for weighted residual empirical processes, resulting in asymptotically distribution-free tests. We investigate the asymptotic properties of this martingale transformation under both the null and alternative hypotheses.
Surprisingly, the martingale-transformation-based test detects local alternatives distinct from the null only at a rate of order $n^{-1/4}$, implying that tests based on the weighted residual empirical process may theoretically outperform their martingale-transformed counterparts in terms of power.
This finding contrasts with existing asymptotically distribution-free tests using martingale-transformed residual-marked empirical processes, which typically achieve a sensitivity rate of order $n^{-1/2}$. Furthermore, we note that when testing for membership in a location class of distributions, the martingale-transformed residual empirical process can also achieve the sensitivity rate of order $n^{-1/2}$ \citep{khmaladze2009}.
On the other hand, this $n^{-1/4}$ rate, while slower, is dimension-independent, unlike local smoothing tests whose sensitivity to local alternatives is often dimension-dependent (e.g., \cite{hardle1993}; \cite{zheng1996}). This dimension-independence is a notable theoretical and practical advantage in high-dimensional settings.
Numerical studies in Section~6 confirm these theoretical findings, showing that the martingale transformation-based test generally has lower power than the corresponding test without this transformation.

%Numerical studies presented in Section 6 validate these theoretical results, showing that tests based on weighted residual empirical processes generally exhibit higher power than the corresponding martingale transformation-based tests. Nevertheless, due to its asymptotically distribution-free property and its ability to mitigate dimensionality-related challenges, this methodology remains highly practical, particularly in scenarios with large sample sizes and high-dimensional parameter and predictor vectors.

The proposed methodology can also be applicable for checking whether the conditional variance function $\sigma^2(\cdot)$ in~(\ref{1.1}) belongs to a given parametric family of functions $\tilde{\mathcal{M}}=\{ \sigma^2(\cdot, \theta): \theta \in \tilde{\Theta} \subset \mathbb{R}^q\}$, when the mean function is in the class $\mathcal{M}=\{ m(\cdot, \beta): \beta \in \Theta \subset \mathbb{R}^p\}$. The asymptotic properties of the resulting weighted residual empirical process are investigated under both the null and alternative hypotheses provided that $ p \vee q=o(n^{1/3} / \log^2{n})$, where $ p \vee q = \max\{p, q\}$.
We also show that the corresponding test without martingale transformation cannot be asymptotically distribution-free in theory, whereas the test with martingale transformation can only be sensitive to local alternatives at the rate of order $n^{-1/4}$.

The rest of this paper is organized as follows. In Sections~2 and 3, we define the weighted residual empirical process and its martingale transformation, and investigate their asymptotic properties under the null hypothesis. Section 4 presents power analyses under the global and local alternative hypotheses. Section 5 discusses the selection of the weight function for the test statistic in practice. Section 6 includes simulation studies and a real data analysis to evidence the performance of the proposed tests. Section 7 contains some concluding remarks and topics for future research. The proofs of all theoretical results in the main text and tests for checking conditional variance functions are deferred to the Supplementary Material for saving space.

\section{Weighted residual empirical processes}
First, we consider the model specification test for the mean function and assume temporarily that the error term $\varepsilon$ is independent of $X$ in model~(\ref{1.1}). The extensions to heteroscedastic cases and dependent data will be discussed in the Supplementary Material.
Our interest here is to check whether the regression function $m(x)=E(Y|X=x)$ belongs to some given parametric class $\mathcal{M}=\{m(\cdot, \beta): \beta \in \Theta \subset \mathbb{R}^p \}$. Therefore, the null and alternative hypotheses are
\begin{eqnarray*}
&H_0: & \mathbb{P}\{m(X) = m(X, \beta_0)\} =1, \quad {\rm for \ some} \ \beta_0 \in \Theta \subset \mathbb{R}^p; \\
&H_1: &\mathbb{P}\{ m(X) \neq m(X, \beta) \} >0, \quad {\rm for \ all} \ \beta \in \Theta \subset \mathbb{R}^p.
\end{eqnarray*}
To illustrate our method, we rewrite the parametric mean function, incorporating a constant intercept term, within the class $\mathcal{M} = \{ m(\cdot, \beta)= \gamma + m_0(\cdot, \vartheta): \beta = (\gamma, \vartheta) \in \mathbb{R}^p\}$.
Here the function $m_0(\cdot, \vartheta)$ does not include the intercept term. In the case with $\gamma=0$, we still treat $\gamma$ as an unknown parameter needed to be estimated.
Write $e=Y-m(X, \tilde{\beta}_0)$, where $\tilde{\beta}_0$ is defined as
\begin{equation}\label{2.1}
\tilde{\beta}_0 = \mathop{\rm arg\ min}\limits_{\beta \in \Theta} E\{ Y-m(X, \beta) \}^2 =\mathop{\rm arg\ min}\limits_{\beta \in \Theta} E\{ m(X)-m(X, \beta) \}^2.
\end{equation}
Let $g(X)$ be a given weight function of $X$ and $g_0(X)=g(X)-E[ g(X)]$. Possible selection of $g$ will be discussed later.
Our methodology is based on the following result which will be proved in Section S1 of the Supplementary Material.

\begin{prop}\label{equ-H0}
Suppose that $\varepsilon$ is independent of $X$ and Condition (A5) in the following holds. Then there exists a weight function $g(\cdot)$ such that
$$ H_0 \ {\rm holds} \ \Longleftrightarrow \ E[g_0 (X)I (e \leq t)] = 0 \ \  {\rm for \ all} \ t \in \mathbb{R}. $$
\end{prop}

Suppose we have an i.i.d. sample $(X_i, Y_i), 1 \leq i \leq n$, following the same distribution with $(X, Y)$.
Motivated by Proposition~\ref{equ-H0}, we define a weighted residual empirical process for testing $H_0$ as
\begin{equation}\label{2.2}
\hat{U}_n(t) = \frac{1}{\sqrt{n}} \sum_{i=1}^{n} [g(X_i)-\bar{g}] I(\hat{e}_i \leq t),
\end{equation}
where $\bar{g}=n^{-1} \sum_{i=1}^{n}g(X_i)$, and $\hat{e}_i= Y_i-m(X_i, \hat{\beta})$ with $\hat{\beta}$ being a $L_2$-norm consistent estimator of $\tilde{\beta}_0$.
In this paper we restrict ourselves to the least squares estimator of $\tilde{\beta}_0$, that is,
\begin{equation}\label{2.3}
\hat{\beta} = \mathop{\rm arg\ min}\limits_{\beta \in \Theta} \sum_{i=1}^n[Y_i-m(X_i, \beta)]^2.
\end{equation}

%The residual empirical processes in parametric regression models have been well studied by \cite{loynes1980}, \cite{portnoy1986}, \cite{mammen1996}, and \cite{koul2002}, among many others. \cite{akritas2001} extended these classical results to the nonparametric regression settings. Another work related to ours is \cite{keilegom2008} which proposed a test for the mean function based on the difference between the parametric and nonparametric residual empirical processes. \cite{escanciano2018} extended this method for testing semiparametric hypotheses in regression models with dependent data and introduced a novel integral transformation that led their test to be asymptotically distribution-free. Note that these two tests require the nonparametric estimation for the unknown mean function to construct the residual empirical process. Thus their methods are different from ours in this paper.

\begin{remark}
We provide a more detailed rationale for using the weighted residual empirical process $\hat{U}_n(t)$ instead of the residual-marked empirical process to address the challenge posed by high dimensionality. Recall that the residual-marked empirical process-based test proposed by \cite{stute1997} is based on $S_n(x) = n^{-1/2} \sum_{i=1}^n \hat{e}_i I(X_i \leq x) $, where the indicator function $I(X_i \leq x)=\prod_{j=1}^d I(X_{ij} \leq x_j)$, and $X_{ij}$ and $x_j$ are the $j$-th component of $X_i$ and $x$, respectively.
A critical issue arises in high-dimensional spaces: $I(X_i \leq x)=0$ occurs whenever at least one component satisfies $I(X_{ij} \leq x_j)=0.$
Due to data sparsity in high dimensions, this condition is highly likely to be met, as noted by \cite{escanciano2006a}. Consequently, the test of \cite{stute1997} based on the residual-marked empirical process suffers severely from the curse of dimensionality.
In contrast, tests based on the weighted residual empirical process $\hat{U}_n(t)$ rely solely on the one-dimensional indicator function  $I(\hat{e}_i \leq t)$, which avoids the data sparsity issue in high-dimensional spaces. By incorporating a suitable weight function $g(\cdot)$, the weighted residual empirical process offers a mechanism to mitigate the dimensionality problem. However, selecting an appropriate weight function is nontrivial, and an unsuitable choice can lead to power loss under alternative hypotheses. Section~5 provides a detailed discussion on how to select a weight function to ensure that the proposed tests maintain good power performance.
\end{remark}

\subsection{Limiting null distributions of weighted residual empirical processes}
In this subsection, we investigate the asymptotic properties of $\hat{U}_n(t)$ under $H_0$, assuming temporarily that the weight function $g(\cdot)$ is given without any unknowns. To this end, we introduce some notations and regularity conditions.
Write $\dot{m}(x, \beta)= \frac{\partial m(x, \beta)}{\partial \beta}$, $\ddot{m}(x, \beta)= \frac{\partial \dot{m}(x, \beta)}{\partial \beta^{\top}}$, and
$\Sigma= E\{\dot{m}(X, \tilde{\beta}_0) \dot{m}(X, \tilde{\beta}_0)^{\top}- [m(X)-m(X, \tilde{\beta}_0)]\ddot{m}(X, \tilde{\beta}_0)\}$.
Throughout this paper, $F(x)$ denotes a measurable envelope function satisfying $E|F(X)|^k = O(1)$ for some integer $k \geq 4$, which will be specified later in Conditions (A9) and (A10). $C$ represents a generic constant independent of $n$, and its value may change across appearances. The following conditions apply for all $x$ in its support domain and all limits are taken as $n \to \infty$. To facilitate the theoretical development, define $F(x) = \max\{1, |f(x)|, |g_0(x)|, |m(x)|, \sigma(x) \}$. Note that $F(\cdot)$ serves as an envelope function for $f(\cdot)$, and the specific function $f(\cdot)$ may vary across different conditions.

(A1) (i) $E(\varepsilon^2) < \infty$; (ii) $|\dot{m}_i(x, \beta)| \leq F(x)$ and $|\ddot{m}_{ij}(x, \beta)| \leq F(x)$ for all $\beta \in U(\tilde{\beta}_0)$, where $U(\tilde{\beta}_0)$ is some neighborhood of $\tilde{\beta}_0$, $\dot{m}_i(x, \beta)$ is the $i$-component of $\dot{m}(x, \beta)$, and $\ddot{m}_{ij}(x, \beta)$ is the $(i,j)$-element of the matrix $\ddot{m}(x, \beta)$.

(A2) Let $\psi(x, \beta)=[m(x)-m(x, \beta)]\dot{m}(x, \beta)$ and $ \dot{\psi}(x, \beta) = \frac{\partial \psi}{\partial \beta^{\top}}(x, \beta)$. (i) $E|\psi_i(X, \tilde{\beta}_0)|^4 \leq C$; (ii) $ |\dot{\psi}_{ij}(x, \beta)| \leq F(x)$ for all $\beta \in U(\tilde{\beta}_0)$, where $\psi_i(x, \tilde{\beta}_0)$ is the $i$-component of $\psi(x, \tilde{\beta}_0)$, and $\dot{\psi}_{ij}(x, \beta)$ is the $(i,j)$-element of the matrix $\dot{\psi}(x, \beta)$.

(A3) Let $\Sigma(\beta) = -E[\dot{\psi}(X, \beta)]$ and let $\Sigma_i(\beta)$ be the $i$-row of $\Sigma(\beta)$. (i) For all $\beta \in U(\tilde{\beta}_0)$, the matrix $\Sigma(\beta)$ is nonsingular and satisfies $0 < C^{-1} \leq \lambda_{\min}\{\Sigma(\beta)\} \leq \lambda_{\max}\{\Sigma(\beta)\} \leq C < \infty$, where $\lambda_{\min}\{\Sigma(\beta)\}$ and $\lambda_{\max}\{\Sigma(\beta)\}$ denote the smallest and largest eigenvalues of $\Sigma(\beta)$, respectively. (ii) The matrix $\dot{\Sigma}_i(\beta) = \frac{\partial \Sigma_i(\beta)}{\partial \beta}$ satisfies $ \max_{1\leq i,l \leq p} |\lambda_l\{\dot{\Sigma}_i(\beta)\}| \leq C $ for all $ \beta  \in U(\tilde{\beta}_0)$,
where $\{\lambda_l\{\dot{\Sigma}_i(\beta)\}: 1 \leq l \leq p\}$ are the eigenvalues of the matrix $\dot{\Sigma}_i(\beta)$.

(A4) For all $\beta_1$ and $\beta_2$ in the neighborhood $ U(\tilde{\beta}_0) $ of $\tilde{\beta}_0$,
\begin{eqnarray*}
&& (\mathrm{i}) \ |\dot{\psi}_{ij}(x, \beta_1)-\dot{\psi}_{ij}(x, \beta_2)| \leq \sqrt{p} \|\beta_1-\beta_2\|F(x), \\
&& (\mathrm{ii}) \ |\ddot{m}_{ij}(x, \beta_1)- \ddot{m}_{ij}(x, \beta_2)| \leq \sqrt{p} \|\beta_1-\beta_2\|F(x),
\end{eqnarray*}
where $\|\cdot\|$ denotes the $L_2$-norm.

(A5) The vector $\tilde{\beta}_0$ lies in the interior of the compact subset $\Theta$ and is the unique minimizer of~(\ref{2.1}).

%\begin{prop}\label{prop1}
%Assume Conditions (A1)-(A5). If $p^2/n \to 0$, then $\|\hat{\beta}-\tilde{\beta}_0\|=O_p(\sqrt{p/n})$. Further, if $(p^3\log^2{n})/n \to 0$, then
%\begin{equation*}
%\sqrt{n} (\hat{\beta}-\tilde{\beta}_0)
%= \frac{1}{\sqrt{n}} \sum_{i=1}^{n}e_i \Sigma^{-1} \dot{m}(X_i, \tilde{\beta}_0)+O_p(\sqrt{\frac{p^3}{n}})
%  +o_p(\frac{p^3 \log{n}}{n}),
%\end{equation*}
%where $e_i = Y_i-m(X_i, \tilde{\beta}_0)$.
%\end{prop}

%Here, the remainders $r_n = o_p(1) \ {\rm or} \ O_p(1)$ means $\|r_n\| = o_p(1) \ {\rm or} \ O_p(1)$, respectively.
%To derive the asymptotic properties of the process $\hat{U}_n(t)$ in diverging dimension settings, we need some additional regularity conditions.

(A6) Assume that
\begin{eqnarray*}
&&{\rm(i)}  \ \lambda_i (E[F^2(X) \dot{m}(X, \tilde{\beta}_0)\dot{m}(X, \tilde{\beta}_0)^{\top}] )  \leq C, \\
&&{\rm(ii)} \ |\lambda_i (E[\ddot{m}(X, \tilde{\beta}_0)]) | \leq C, \quad  |\lambda_i (E[\dot{m}_j(X, \tilde{\beta}_0)\ddot{m}(X,
              \tilde{\beta}_0)]) | \leq C, \\
&&{\rm(iii)}\ |\lambda_i (E[g_0(X)f_{\varepsilon}(t + m(X, \tilde{\beta}_0)-m(X)) \ddot{m}(X, \tilde{\beta}_0)])| \leq C \ \ \forall \ t,
\end{eqnarray*}
where $\{\lambda_i(M)\}_{i=1}^p$ denotes the eigenvalues of a $p \times p$ matrix $M$ and $f_{\varepsilon}(\cdot)$ is the density function of $\varepsilon$.

(A7) The weight function $g(x)$ satisfies that $g(X)$ is sub-exponential with a parameter $\sigma_0 >0$, that is, $E[\exp\{t (g(X)-Eg(X)) \}] \leq \exp(\sigma_0^2 t^2 /2)$ for all $|t| < 1/\sigma_0$.

Conditions (A1)-(A5) are primarily used to derive the $L_2$-norm consistency and asymptotically linear expansion of $\hat{\beta} - \tilde{\beta}_0$. To save space, we present these asymptotic properties of $\hat{\beta} - \tilde{\beta}_0$ in the Supplementary Material (see Propositions S1 and S2).
\cite{tan2022} previously established the consistency and asymptotically linear expansion of $\hat{\beta} - \tilde{\beta}_0$ in diverging-dimension scenarios. Their work required the existence of the fourth moment of the response and the third derivative of the regression function $m(x, \beta)$. Our research significantly relaxes these assumptions, requiring only the second moment of the error term and the second derivative of $m(x, \beta)$.
Condition (A6) is similar to the regularity condition (B1) in \cite{tan2022}, which is used to control the convergence rate of remainders in diverging dimensional statistical inference, see also \cite{fanpeng2004} and \cite{zhang2008} for instance. Condition (A7) is a technical condition needed for controlling the divergence rate of the dimension $p$. If this condition is violated, we may require a slower divergence rate of $p$ such as $p = o(n^{1/(3+c)})$, to ensure the convergence of related high-dimensional empirical processes involved in the decomposition of $\hat{U}_n(t)$. Here $c$ is some small positive constant.
In the case of testing linear regression models, i.e., $m(x, \beta) =\beta^{\top}x$, Conditions (A4), (A6)(ii), and (A6)(iii) hold automatically because $\dot{m}(x, \beta) = x$ and $\ddot{m}(x, \beta) = 0$. These conditions are necessary to control the convergence rate of the remainders arising from general nonlinear regression models in the proofs of related empirical processes. Further examples illustrating these conditions are provided in Section S5 of the Supplementary Material.

The following result establishes the asymptotic properties of $\hat{U}_n(t)$ under both the null and alternative hypotheses. Let $F_{\varepsilon}(t)$ denote the cumulative distribution function of $\varepsilon$ and define $\nu_n = \max\{1, |\varphi_e(-c_{1n})|, |\varphi_e(c_{2n})|\}$, where $\varphi_e(t) = \dot{f}_e(t)/f_e(t)$ with $f_e(t)$ being the density function of $e$ and $\dot{f}_{e}(t) = \frac{df_{e}(t)}{dt}$. Here, $c_{1n}, c_{2n}$ $ \uparrow \infty$ are the truncated parameters that will be introduced in Section 3.

\begin{theorem}\label{expan-Un}
Assume Conditions (A1)-(A7). If $p^3\log^6{n} = o(n)$, then, uniformly in $t \in \mathbb{R}$,
\begin{eqnarray}\label{2.4}
\hat{U}_n(t)
&=& \frac{1}{\sqrt{n}} \sum_{i=1}^n g_0(X_i) \{I[\varepsilon_i \leq t + m(X_i, \tilde{\beta}_0)-m(X_i)]
    - E[F_{\varepsilon}(t+m(X, \tilde{\beta}_0)-m(X))]\} \nonumber \\
&&  + \frac{1}{\sqrt{n}} \sum_{i=1}^{n}e_i M(t)^{\top} \Sigma^{-1} \dot{m}(X_i, \tilde{\beta}_0)+R_n(t),
\end{eqnarray}
where $M(t)=E[g_0(X)\dot{m}(X, \tilde{\beta}_0) f_{\varepsilon}(t+m(X, \tilde{\beta}_0)-m(X))]$, and $R_n(t)$ is a remainder such that $\sup_{t \in \mathbb{R}}|R_n(t)|=o_p(1)$.
Furthermore, if $p^3 (\nu_n\log{n})^6 = o(n)$, then (\ref{2.4}) continues to hold with the remainder satisfying $\sup_{t \in \mathbb{R}}|R_n(t)|=o_p(\nu_n^{-1})$.
\end{theorem}

We acknowledge that \cite{stute2008} established a decomposition of the weighted residual empirical process $\hat{U}_n(t)$ in fixed-dimension scenarios. Here we extend this decomposition to diverging-dimension settings.
The rate $\sup_t|R_n(t)|=o_p(\nu_n^{-1})$ is used to ensure $\hat{T}_n R_n(t) = o_p(1)$ uniformly in $t$ where $\hat{T}_n$ is the martingale transformation proposed in Section 3.

Recall that under $H_0$, we have $\tilde{\beta}_0 = \beta_0$, $m(X) = m(X, \tilde{\beta}_0)$, and $\Sigma= E[\dot{m}(X, \beta_0) \dot{m}(X, \beta_0)^{\top}]$. It follows from Theorem \ref{expan-Un} that under $H_0$,
\begin{eqnarray}\label{2.6}
\hat{U}_n(t)
&=&  \frac{1}{\sqrt{n}} \sum_{i=1}^{n} g_0(X_i) [I(\varepsilon_i \leq t)-F_{\varepsilon}(t)] + f_{\varepsilon}(t)
     \frac{1}{\sqrt{n}} \sum_{i=1}^{n} \varepsilon_i M^{\top}\Sigma^{-1}\dot{m}(X_i, \beta_0)+ R_n(t) \nonumber \\
&=:& U_n^1(t) + R_n(t),
\end{eqnarray}
where $M=E[g_0(X)\dot{m}(X, \beta_0)]$ and $\sup_t|R_n(t)| = o_p(1)$.
When $H_0$ is false, we can choose a suitable weight function $g(X)$ such that the limit of $E[g_0(X) F_{\varepsilon}(t+m(X, \tilde{\beta}_0)-m(X))]$ is nonzero. Consequently, the first sum in~(\ref{2.4}) may diverge to infinity in probability as $n \to \infty$. Since the second sum in~(\ref{2.4}) is bounded in probability uniformly in $t$, it follows that tests based on the process $\hat{U}_n(t)$ could be consistent with asymptotic power $1$ under $H_1$. Therefore, we can construct a test for $H_0$ based on a functional of $\hat{U}_n(t)$.
In this paper, we consider the Cra\'{m}er-von Mises type test with a properly chosen weight function $g$ based on $\hat{U}_n(t)$ as
\begin{eqnarray*}
CvM_n = \int_{\mathbb{R}} |\hat{\rho}^{-1} \hat{U}_n(t)|^2 d\hat{F}_{\hat{e}}(t),
\end{eqnarray*}
where $\hat{\rho}^2 = n^{-1}\sum_{i=1}^n [g(X_i)-\bar{g}]^2$ and $\hat{F}_{\hat{e}}(t)$ denotes the empirical distribution function of $\{ \hat{e}_i \}_{i=1}^n$.
The next theorem provides the limiting null distribution of $\hat{U}_n(t)$ in diverging-dimension scenarios. Its proof will be provided in the Supplementary Material.

\begin{theorem}\label{hatU-limitnull}
Assume Conditions (A1)-(A7). If $p^3 \log^6{n} = o(n)$, then under $H_0$,
\begin{eqnarray*}
\hat{U}_n(t) \longrightarrow U_{\infty}^1(t) \quad { in \ distribution},
\end{eqnarray*}
in the Skorohod space $D[-\infty, \infty].$ Here $U_{\infty}^1(t)$ is a zero-mean Gaussian process with the covariance function $K(s, t)$ which is the pointwise limit of $K_n(s, t)$ given by
\begin{eqnarray*}
&&  K_n(s, t) \\
&=& Cov(U_n^1(s), U_n^1(t)) \\
&=& \rho_n^2 [F_{\varepsilon}(s \wedge t)-F_{\varepsilon}(s)F_{\varepsilon}(t)] + M^{\top} \Sigma^{-1} M
    \{f_{\varepsilon}(s) E[\varepsilon I(\varepsilon \leq t)] + f_{\varepsilon}(t) E[\varepsilon I(\varepsilon \leq s)] + f_{\varepsilon}(s)f_{\varepsilon}(t) E(\varepsilon^2) \},
\end{eqnarray*}
where $s \wedge t = \min\{s,t\}$ and $\rho_n^2 = E[g_0(X)]^2$.
\end{theorem}

Note that the covariance function $K_n(s, t)$ involves the second moment of the error term $ E(\varepsilon^2)$. Therefore, it is necessary to impose the existence of the second moment of the error term in Condition (A1) to derive the Gaussian process approximation of $\hat{U}_n(t)$.
Theorem \ref{hatU-limitnull} and  the Extended Continuous Mapping Theorem (e.g., Theorem 1.11.1 of \cite{van1996}) yield the limiting null distribution of $CvM_n$.

\begin{theorem}\label{CVM-limitnull}
Assume Conditions (A1)-(A7). If $p^3 \log^6{n} = o(n)$, then under $H_0$,
\begin{eqnarray*}
CvM_n \longrightarrow \int_{\mathbb{R}} |\rho^{-1} U_{\infty}^1(t)|^2 dF_{\varepsilon}(t) \quad in \ distribution,
\end{eqnarray*}
where $\rho$ is the limit of $\rho_n$ and $U_{\infty}^1(t)$ is given in Theorem \ref{hatU-limitnull}.
\end{theorem}

\subsection{Smooth residual bootstrap approximation}
According to Theorems \ref{hatU-limitnull} and \ref{CVM-limitnull}, it is readily seen that the process
$U_{\infty}^1(t)$ involves the unknown functions $F_{\varepsilon}(t)$ and $f_{\varepsilon}(t)$,
which make the test $CvM_n$ infeasible for critical value determination in practice. Therefore, we suggest a smooth residual bootstrap approximation for the limiting null distribution of $CvM_n$.
This bootstrap method was first proposed by \cite{koul1994} to approximate the residual empirical process for linear models in fixed-dimension scenarios. \cite{neumeyer2009} extended this method to nonparametric regression models, while \cite{mora2005} and \cite{dette2007} adopted this approach to approximate the limiting null distributions of their tests. We now extend this approach to handle diverging-dimension scenarios. The algorithm is as follows:

1. Generate the bootstrap errors as $\varepsilon_i^{\ast} = \tilde{\varepsilon}_i^* + v_n Z_i$, where $\tilde{\varepsilon}_i^*$ are randomly sampling with replacement from the centered residuals $\{ \tilde{e}_i = \hat{e}_i- n^{-1} \sum_{i=1}^n \hat{e}_i:  i=1, \cdots, n \}$,
$v_n$ is a smoothing parameter, and $Z_1, \cdots, Z_n$ are independent, centered random variables with density $l(\cdot)$, independent of the original sample $\mathcal{Y}_n=\{(X_1, Y_1), \cdots, (X_n, Y_n)\}$.

2. Generate a bootstrap sample according to the model $ Y_i^{\ast}=m(X_i, \hat{\beta}) + \varepsilon_i^{\ast}, i=1, \cdots, n$. Let $\hat{\beta}^*$ be the bootstrap estimator obtained by the least squares method based on the bootstrap sample $\{(X_j, Y_j^*)\}_{j=1}^n$.

3. Define the bootstrap version of the test statistic as
$$ CvM_n^* = \int_{\mathbb{R}} |\hat{\rho}^{-1} \hat{U}_n^*(t)|^2 d\hat{F}_{\hat{\varepsilon}^*}(t), $$
where %$\hat{\rho^*}^2 = n^{-1}\sum_{i=1}^n [\hat{g}^*(X_i)-\bar{\hat{g}}^*]^2$,
$ \hat{U}_n^*(t) = n^{-1/2}\sum_{i=1}^n [g(X_i)-\bar{g}] I(\hat{\varepsilon}_i^{\ast} \leq t)$ with $\hat{\varepsilon}_i^{\ast}=Y^{*}_i-m(X_i, \hat{\beta}^*)$,
%$\hat{g}^*(\cdot)$ is defined in the same way as $\hat{g}_l(\cdot)$ in (\ref{6.5}) except the estimators being replaced by the corresponding bootstrap versions,
and $\hat{F}_{\hat{\varepsilon}^*}(t)$ is the empirical distribution of $\{\hat{\varepsilon}_i^{\ast}\}_{i=1}^n$.

4. Repeat steps 1-3 a large number of times, say $B$ times. For a given nominal level $\zeta$, the critical value is determined by the upper $\zeta$-quantile of the bootstrap distribution of $CvM_n^*$.

To establish the validity of the smooth residual bootstrap, we need some additional regularity conditions.

(A8) (i) The kernel density function $l(\cdot)$ in Step 1 is positive, symmetric, and twice continuously differentiable, satisfying $\int_{\mathbb{R}} tl(t) dt =0$ and $\int_{\mathbb{R}} t^4 l(t) dt < \infty$. (ii) The smoothing parameter $v_n$ satisfies $v_n = o(1)$ and $\log{n} = o(nv_n^4)$.

These conditions are similar to Assumptions 4-5 as presented in \cite{mora2005}.
In the simulation studies, the chosen sample $\{Z_1, \cdots, Z_n\}$ in Step 1 consists of i.i.d. standard normal random variables, which ensures that the density function $l(\cdot)$ satisfies Condition (A8).
The following theorem establishes the asymptotic validity of the smooth residual bootstrap approximation under $H_0$ in diverging-dimension settings. Its proof is provided in the Supplementary Material. Define $c_{n, \zeta}^* $ by $\mathbb{P}(CvM_n^* > c_{n, \zeta}^* | \mathcal{Y}_n) = \zeta$.

\begin{theorem}\label{bootstrap-null}
Assume that Conditions (A1)-(A8) hold. If $p^3\log^6{n}=o(n)$, then, under $H_0$ and conditionally on the original sample $\mathcal{Y}_n$,
$$ CvM_n^* \longrightarrow \int_{\mathbb{R}} |\rho^{-1} U_{\infty}^{1*}(t)|^2 dF_{\varepsilon}(t) \quad {in \ distribution} $$
in probability, where $U_{\infty}^{1*}(t)$ has the same distribution as the Gaussian process $U_{\infty}^1(t)$ given in Theorem \ref{hatU-limitnull}. This implies that
$ \lim_{n \to \infty} \mathbb{P}( CvM_n > c_{n, \zeta}^*) = \zeta. $
\end{theorem}

\section{Martingale transformation}
%It is well known that tests based on functionals of $\hat{U}_n(t)$, such as the Kolmogorov-Smirnov type test and the Cra\'{m}er-von Mises type test, cannot be asymptotically distribution-free. Therefore, we need to resort to resampling methods, such as the wild bootstrap, to approximate the limiting null distribution of such tests (e.g., \cite*{stute1998a}).
When the dimension $p$ is large, the computation of the bootstrap procedure becomes cumbersome and time-consuming, as demonstrated by the computational time of $CvM_n$ based on the residual smoothing bootstrap in our simulation studies. Therefore, an asymptotically distribution-free test is desirable in high dimension settings.
For this purpose, we examine the asymptotic expansion of the process $\hat{U}_n(t)$ under $H_0$ more carefully.
The leading term $U^1_n(t)$ of $\hat{U}_n(t)$ in (\ref{2.6}) can be decomposed as
$$U_n^1(t) = U_n^0(t) - F_{\varepsilon}(t) \frac{1}{\sqrt{n}} \sum_{i=1}^{n} g_0(X_i) + f_{\varepsilon}(t) \frac{1}{\sqrt{n}} \sum_{i=1}^{n} \varepsilon_i M^{\top}\Sigma^{-1}\dot{m}(X_i, \beta_0), $$
where $U_n^0(t) = n^{-1/2} \sum_{i=1}^{n} g_0(X_i) I(\varepsilon_i \leq t)$.
It is readily seen that $ cov\{U_n^0(s), U_n^0(t)\} = \rho_n^2 F_{\varepsilon}(s \wedge t) $.
Consequently, $U_n^0(t)$ admits the same covariance structure as a scaled time-transformed Brownian motion $ \rho_n B(F_{\varepsilon}(t))$, where $B(t)$ is the standard Brownian motion. The convergence of the finite-dimensional distribution and the asymptotic tightness of $U_n^0(t)$ can be proved by the standard arguments. Then it yields
$$ U_n^0(t) \longrightarrow \rho B(F_{\varepsilon}(t)) \quad {\rm in \ distribution},$$
in the Skorohod space $D[-\infty, \infty]$, where $\rho$ is the limit of $\rho_n$.

Note that functionals of $U_n^0(t)$, such as the Cram\'{e}r-von Mises functional and the Kolmogorov-Smirnov functional, after removing the parameter $\rho$ and using the time transformation $z=F_{\varepsilon}(t)$, are asymptotically pivotal and have tractable limiting distributions.
Therefore, in order to develop an asymptotically distribution-free test based on $\hat{U}_n(t)$, we transform $\hat{U}_n(t)$ to another process that would admit the same limiting distribution as $U_n^0(t)$ and simultaneously eliminate the shift terms in $\hat{U}_n(t)$.
Inspired by the martingale transformation introduced by \cite*{stute1998b} for model checking, a martingale transformation in our settings can also achieve this goal.
This methodology can be traced back to \cite{khmaladze1981} for testing the goodness-of-fit of cumulative distribution functions.
The basic principle of martingale transformation is as follows.

Let $W_{\varepsilon}(t) = (F_{\varepsilon}(t), f_{\varepsilon}(t))^{\top}$ and let $h_{\varepsilon}(t)=\frac{\partial W_{\varepsilon}(t)}{\partial F_{\varepsilon}(t)} $ be the Radon-Nikodym derivative of $W_{\varepsilon}(t)$ with respect to $F_{\varepsilon}(t)$. It is easy to see that $ h_{\varepsilon}(t)=(1, \varphi_{\varepsilon}(t))^{\top}$ with $\varphi_{\varepsilon}(t) = \dot{f}_{\varepsilon}(t)/f_{\varepsilon}(t)$. Furthermore, consider the matrix
$ \Gamma_{\varepsilon}(t)=\int_{t}^{\infty} h_{\varepsilon}(z)h_{\varepsilon}(z)^{\top} dF_{\varepsilon}(z)$.
The martingale transformation is defined as
\begin{equation}\label{3.1}
T_{\varepsilon}f(t)=f(t)-\int_{-\infty}^t h_{\varepsilon}(z)^{\top} \Gamma_{\varepsilon}(z)^{-1} \int_z^{\infty} h_{\varepsilon}(v) df(v) dF_{\varepsilon}(z)
\end{equation}
where $f(t)$ is either a bounded variation function or a stochastic process such that the integral in~(\ref{3.1}) is well defined. It may also be a vector of functions sometimes.
Usually the matrix $\Gamma_{\varepsilon}(t)$ is assumed to be nonsingular for $t \in \mathbb{R}$ (e.g., \cite*{stute1998b}, \cite{bai2003}, and \cite{tan2019a}). Although the matrix $\Gamma_{\varepsilon}(t)$ is indeed nonsingular for most of density functions $f_{\varepsilon}(t)$, there are some densities that do not to satisfy this assumption; see \cite{khmaladze2009} for some examples. In such cases, $\Gamma_{\varepsilon}(t)^{-1}$ is the generalized inverse of $\Gamma_{\varepsilon}(t)$ satisfying the first Moore-Penrose assumption $\Gamma_{\varepsilon}(t)\Gamma_{\varepsilon}(t)^{-1}\Gamma_{\varepsilon}(t) = \Gamma_{\varepsilon}(t). $
Interestingly, \cite{nikabadze1997} showed that the transformation $T_{\varepsilon}$ is unique irrespective of whether the generalized inverse of $\Gamma_{\varepsilon}(t)$ is used.

Recall that our purpose is to replace the process $\hat{U}_n(t)$ by the transformed process which shares the same limit distribution as $U_n^0(t)$ and ensures the shift terms in $\hat{U}_n(t)$ vanish.
The next result shows that the martingale transformation in (\ref{3.1}) can achieve these two goals simultaneously. Its proof will be given in the Supplementary Material.

\begin{prop}\label{martin-prop}
Assume that $\Gamma_{\varepsilon}(z)$ is nonsingular for any $z \in \mathbb{R}$. Then the martingale transformation $T_{\varepsilon}$ satisfies the following two properties: \\
{\rm (i)}  $T_{\varepsilon} F_{\varepsilon}(t) \equiv 0$ and $T_{\varepsilon} f_{\varepsilon}(t) \equiv 0$; \\
{\rm (ii)} $T_{\varepsilon}U_n^0(t)$ and $U_n^0(t)$ admit the same weak limit in distribution.
\end{prop}

Since the martingale transformation $T_{\varepsilon}$ is a linear operator, Proposition \ref{martin-prop} implies that
\begin{equation}\label{3.3}
T_{\varepsilon}U_n^1(t)=T_{\varepsilon}U_n^0(t) \longrightarrow \rho B(F_{\varepsilon}(t)) \quad {\rm in \ distribution}
\end{equation}
in the Skorohod Space $D[-\infty, \infty)$. Consequently, we get rid of the shift terms in $U_n^1(t)$, and the transformed process $T_{\varepsilon} U_n^1(t)$ weakly converges to a technically simple process. Furthermore, under some regularity conditions and $p^3 (\nu_n\log{n})^6 = o(n)$, we can show that the transformed remainder process $T_{\varepsilon} R_n(t)= o_p(1)$ uniformly in $t$.
Therefore, the martingale transformed process $T_{\varepsilon} \hat{U}_n(t)$, after dropping off the parameter $\rho$ and using the transformation $z=F_{\varepsilon}(t)$, weakly converges to the standard Brownian motion $B(t)$.

The above transformation $T_{\varepsilon}$ is still not applicable for statistical inference as it involves some unknowns such as $f_{\varepsilon}(t)$ and $F_{\varepsilon}(t)$. We should replace these functions with corresponding estimators to get a final process $\hat{T}_n \hat{U}_n(t)$. As we do not make any assumption except for the smoothness of the density function $f_{\varepsilon}(t)$, it should be estimated nonparametrically. For instance, a standard Nadaraya-Watson estimator for $f_{\varepsilon}(t)$ can be adopted:
$$ \hat{f}_{\hat{e}}(t) = \frac{1}{nh} \sum_{i=1}^{n} K(\frac{t-\hat{e}_i}{h}), $$
where $\hat{e}_i= Y_i - m(X_i, \hat{\beta})$, $K(\cdot)$ denotes a univariate kernel function and $h$ is the bandwidth. We then obtain the estimators $\hat{h}_{\hat{e}}(t)$ and $\hat{\Gamma}_{\hat{e}}(t)$ of $h_{\varepsilon}(t)$ and $\Gamma_{\varepsilon}(t)$ respectively:
$$  \hat{h}_{\hat{e}}(t) = (1, \hat{\varphi}_{\hat{e}}(t))^{\top} \quad {\rm and} \quad
\hat{\Gamma}_{\hat{e}}(t) = \int_t^{\infty} \hat{h}_{\hat{e}}(z) \hat{h}_{\hat{e}}(z)^{\top} d \hat{F}_{\hat{e}}(z).
$$
Here $\hat{\varphi}_{\hat{e}}(t) = \hat{\dot{f}}_{\hat{e}}(t)/\hat{f}_{\hat{e}}(t)$, $\hat{\dot{f}}_{\hat{e}}(t)$ is the derivative of $\hat{f}_{\hat{e}}(t)$ with respect to $t$, and $\hat{F}_{\hat{e}}(t)$ denotes the empirical distribution function of $\{ \hat{e}_i \}_{i=1}^n$.
An estimator $\hat{T}_n$ for the martingale transformed process $T_{\varepsilon}$ is subsequently obtained by substituting $\hat{h}_{\hat{e}}(t)$ and $\hat{\Gamma}_{\hat{e}}(t)$ in (\ref{3.1}).

Note that $f_{\varepsilon}(t)$ and $\hat{f}_{\hat{e}}(t)$ are in the denominators of $\varphi_{\varepsilon}(t)$ and $\hat{\varphi}_{\hat{e}}(t)$, respectively. If $\inf_{t \in \mathbb{I}} f_{\varepsilon}(t) >0$ with $\mathbb{I}$ being the support of $f_{\varepsilon}(t)$, then the estimated martingale-transformed process $\hat{T}_n \hat{U}_n(t)$ can be directly used to construct the test statistic. However, if $\inf_{t \in \mathbb{I}} f_{\varepsilon}(t) = 0$, then $\hat{\varphi}_{\hat{e}}(t) = \hat{\dot{f}}_{\hat{e}}(t)/\hat{f}_{\hat{e}}(t)$ may be unbounded on the support $\mathbb{I}$. This could cause some theoretical difficulties in deriving the asymptotic properties of the estimated martingale transformation $\hat{T}_n \hat{U}_n(t)$.
We then propose a truncated version of the martingale transformation to overcome this difficulty.
To illustrate the method of truncation, we assume the support of $f_{\varepsilon}(\cdot)$ is the whole real line $\mathbb{R}$. For other cases of support sets, the martingale transformation with truncation and the corresponding theoretical justification can be done very similarly.
The truncated estimators of $h_{\varepsilon}(t)$ and $\Gamma_{\varepsilon}(t)$ are defined as:
$$ \hat{h}_{\hat{e}}^{tr}(t) = (1, \hat{\varphi}_{\hat{e}}(t))^{\top} I(-c_{1n} \leq t \leq c_{2n}) \quad {\rm and} \quad
\hat{\Gamma}_{\hat{e}}^{tr}(t) = \int_t^{\infty} \hat{h}_{\hat{e}}^{tr}(z) \hat{h}_{\hat{e}}^{tr}(z)^{\top} d\hat{F}_{\hat{e}}(z).
$$
Here, $\hat{\varphi}_{\hat{e}}(t) = \frac{\hat{\dot{f}}_{\hat{e}}(t)}{\hat{f}_{\hat{e}}(t) \vee a_n} $, $a\vee b = \max\{a, b\}$, $ c_{1n}, c_{2n} \uparrow \infty $, and $a_n \downarrow 0$ are positive sequences that will be specified later. Using $a_n$ is to avoid zero values of denominator in $\hat{\varphi}_{\hat{e}}(t)$.
Therefore, a truncated estimator of the martingale transformation $T_{\varepsilon}$ is
\begin{equation}\label{3.4}
\hat{T}_n \hat{U}_n(t)=\hat{U}_n(t) - \int_{-\infty}^t \hat{h}_{\hat{e}}^{tr}(z)^{\top} \hat{\Gamma}_{\hat{e}}^{tr}(z)^{-1} \int_z^{\infty} \hat{h}_{\hat{e}}^{tr}(v) d\hat{U}_n(v) d\hat{F}_{\hat{e}}(z).
\end{equation}
To facilitate the analysis for the truncated process $\hat{T}_n \hat{U}_n(t)$, we introduce additional notations.
Write $ W_{\varepsilon}^{tr}(t) = (F_{\varepsilon}(t), f_{\varepsilon}(t))^{\top} I(-c_{1n} \leq t \leq c_{2n})$,
$h_{\varepsilon}^{tr}(t) = \frac{\partial W_{\varepsilon}^{tr}(t)}{\partial F_{\varepsilon}(t)}$,
and
$ \Gamma_{\varepsilon}^{tr}(t) = \int_{t}^{\infty} h_{\varepsilon}^{tr}(z)h_{\varepsilon}^{tr}(z)^{\top} dF_{\varepsilon}(z). $
Then the truncated martingale transformation is
$$ T^{tr}_{\varepsilon} f(t)=f(t)- \int_{-\infty}^t h^{tr}_{\varepsilon}(z)^{\top} \Gamma^{tr}_{\varepsilon}(z)^{-1} \int_z^{\infty} h^{tr}_{\varepsilon}(v) df(v) dF_{\varepsilon}(z). $$
Similar to the arguments for the martingale transformation $T_{\varepsilon}$ in Proposition \ref{martin-prop}, we can show that $ T^{tr}_{\varepsilon}F_{\varepsilon}(t) = o(1)$, $T^{tr}_{\varepsilon}f_{\varepsilon}(t) = o(1)$, and
\begin{equation*}%\label{3.5}
T^{tr}_{\varepsilon}U_n^0(t) \longrightarrow \rho B(F_{\varepsilon}(t)) \quad {\rm in \ distribution}
\end{equation*}
in the Skorohod Space $D[-\infty, \infty)$. This means that the truncated martingale-transformed process $T^{tr}_{\varepsilon}U_n^0(t)$ shares the same limiting distribution as $U_n^0(t)$ and simultaneously ensures the shift terms in
$ \hat{U}_n(t) $ vanish as $n$ tends to infinity.

The following regularity conditions are needed to derive the asymptotic properties of the truncated process $\hat{T}_n \hat{U}_n(t)$ under $H_0$ in diverging-dimension settings.

(A9) The kernel function $K(\cdot)$ satisfies that, (i) $K(\cdot)$ is continuous on $\mathbb{R}$ and has a continuous derivative in its support $[-1, 1]$;
(ii) $K(x)=K(-x)$; (iii) $K(\cdot)$ and $\dot{K}(\cdot)$ are bounded variation; (iv) $\int_{-1}^1 K(u)du = 1$ and $\int_{-1}^1 u^{i} K(u)du = 0$ for $i=1, \cdots, k-1$.

(A10) The density function $f_e(t)$ admits $k-1$ order derivative in $t$ and let $\dot{f}_e(t) =f^{(1)}_e(t) =\frac{df_e(t)}{dt}$, $\ddot{f}_e(t) = f^{(2)}_e(t) = \frac{d^2f_e(t)}{dt^2}$, and $f^{(i)}_e(t) = \frac{d^i f_e(t)}{dt^i}$ for $i=3, \cdots, k-1$.
The function $f^{(k-1)}_e(t)$ satisfies the Lipschitz condition
$$ |f^{(k-1)}_e(t + u) - f^{(k-1)}_e(t)| \leq L|u|, \quad \forall \ u \in U $$
for some neighborhood $U$ of zero.
Moreover, $\sup_t |f^{(i)}_e(t)| < \infty$ and $\int |\frac{f^{(i)}_{e}(t)}{f_e(t)}|^4 d F_e(t) < \infty $ for $i=1, \cdots, k-1$.

(A11) Write $\tau_n = \sum_{i=0}^1 \int_{-c_{1n}}^{c_{2n}} |\frac{f^{(i)}_{e}(t)}{f_e(t)}|^2 dt $ and $\kappa_n = \inf_{-c_{1n} \leq t \leq c_{2n}} f_{e}(t) $. Assume that $a_n/\kappa_n = o(1)$ and the bandwidth $h$ satisfies that $\tau_n \kappa_n^{-1} \sqrt{n} h^{2(k-1)} = o(1)$ and $ \tau_n \kappa_n^{-1} \frac{\log{n}}{\sqrt{n} h^4} =o(1)$.

(A12) The matrix $\Gamma_e(z)$ satisfies that $\inf_{z\leq z_0}|\det(\Gamma_e(z))| > 0$ for any $z_0 \in \mathbb{R}$, where $\det(\Gamma_e(z))$ denotes the determinant of the matrix $\Gamma_e(z)$.

(A13) For some $\delta$ large enough, the function $\varphi_{e}(t)$ has a bounded variation on $[-\delta, \delta]$, and is monotonic on $(-\infty, -\delta]$ and $[\delta, \infty)$, respectively.

Conditions (A9)-(A10) are usually used in the literature of high order nonparametric estimation; see Chapters 2 and 4 of \cite{rao1983}, \cite{zhu1996}, and Chapter 1 of \cite{li2007}, for instance.
Condition (A11) specifies the divergence rate of the truncated parameters $a_n, c_{1n}$ and $c_{2n}$ that are used to control the convergence rate of the remainders in the decomposition of the truncated martingale-transformed process $\hat{T}_n \hat{U}_n(t)$.
If $c_{1n}$ and $c_{2n}$ diverge to infinity slowly enough, it can be shown that these conditions in (A11) for the bandwidth $h$ can be easily satisfied by some common density functions, such as the standard normal density function and the Laplace density function.
Condition (A12) is necessary for the uniformly boundedness of $ \| \Gamma_{e}(z)^{-1} \|$ from infinity.
Condition (A13) is used to ensure the uniform infinitesimality of the transformed remainder process $\hat{T}_n R_n(t)$ in diverging-dimension scenarios.
%Moreover, we show that the conditions for bandwidth $h$ in (A10) can be easily satisfied by some common density functions, such as the standard normal density function and the Laplace density function. For the simplicity of notation, we assume without loss of generality $c_{1n}=c_{2n}=c_n$. Considering the standard normal density function $f_e(t) = \frac{1}{\sqrt{2 \pi}} \exp(-\frac{t^2}{2})$, it is easy to see that $\kappa_n = \frac{1}{\sqrt{2 \pi}} \exp(-\frac{c_n^2}{2}) $ and $\tau_n = 2 c_n + \frac{2}{3} c_n^3.$
%Set $c_n^2 = \frac{1}{4} \log{n}$ for instance. It follows that $\kappa_n = a_1 n^{-1/8}$ and $\tau_n = a_2\log^{3/2}{n}$ for some constants $a_1$ and $a_2$. This leads to the requirements for the bandwidth $h$ such that $(n^{5/8} \log^{3/2}{n}) h^{2(k-1)} = o(1)$ and $\frac{\log^{5/2}{n}}{n^{3/8} h^4} =o(1)$. Similarly, for the Laplace density function $f_e(t) = \frac{1}{2} \exp(-|t|)$, we have
%$ \kappa_n = \frac{1}{2} \exp(-c_n)$ and $ \tau_n = 4 c_n$.
%If $c_n=\frac{1}{8} \log{n}$, then we have $\kappa_n = \frac{1}{8} n^{-1/8}$ and $\tau_n = \frac{1}{2} \log{n}$. Now the conditions for bandwidth $h$ becomes $ (n^{5/8} \log{n}) h^{2(k-1)} = o(1) $  and $ \frac{\log^2{n}}{n^{3/8} h^4} =o(1)$. Thus, if $c_n$ diverges to infinity slowly enough, then the conditions for $h$ in (A10) can be easily satisfied.

The next theorem derives the limiting null distribution of the truncated martingale-transformed process $\hat{T}_n \hat{U}_n(t)$. Its proof is provided in the Supplementary Material.

\begin{theorem}\label{th3.1}
Assume that Conditions (A1)-(A7) and (A9)-(A13) hold and $\Gamma^{tr}_{\varepsilon}(z)$ is nonsingular for each $z \in \mathbb{R}$. If $ p^3 (\nu_n\log{n})^6 = o(n) $, then under $H_0$,
\begin{equation*}
\hat{T}_n \hat{U}_n(t)  \longrightarrow \rho B(F_{\varepsilon}(t)) \quad { in \ distribution},
\end{equation*}
in the Skorohod Space $D[-\infty, \infty)$.
\end{theorem}

For the Cram\'{e}r-von Mises type test based on the truncated martingale-transformed process $\hat{T}_n \hat{U}_n(t)$, following the idea of \cite*{stute1998b}, the (informal) test statistic based on $\hat{T}_n \hat{U}_n(t)$ is defined as
$$ TCvM_n^0 = \int_{-\infty}^{t_0} |\hat{\rho}^{-1}\hat{T}_n \hat{U}_n(t)|^2 d \hat{F}_{\hat{e}}(t). $$
By Theorem \ref{th3.1} and the Extended Continuous Mapping Theorem (e.g., Theorem 1.11.1 of van der Vaart and Wellner (1996)), we have under $H_0$,
$$ TCvM_n^0 \longrightarrow \int_{-\infty}^{t_0} |B(F_{\varepsilon}(t))|^2 dF_{\varepsilon}(t) \quad  {\rm in \ distribution}, $$
where $B(t)$ is the standard Brownian motion.
Since $B(t F_{\varepsilon}(t_0))/\sqrt{F_{\varepsilon}(t_0)} = B(t)$ in distribution, we have
$ \int_{-\infty}^{t_0} |B(F_{\varepsilon}(t))|^2 dF_{\varepsilon}(t) = F_{\varepsilon}^2(t_0) \int_0^1 \frac{|B(t F_{\varepsilon} (t_0))|^2}{F_{\varepsilon}(t_0)} dt = F_{\varepsilon}^2(t_0) \int_0^1 B(t)^2 dt $
in distribution.
Therefore, the final Cram\'{e}r-von Mises type test statistic based on $\hat{T}_n \hat{U}_n(t)$ is
\begin{eqnarray*}
TCvM_n = \frac{1}{\hat{\rho}^2 \hat{F}_{\hat{e}}(t_0)^2} \int_{-\infty}^{t_0} | \hat{T}_n \hat{U}_n(t)|^2 d \hat{F}_{\hat{e}}(t).
\end{eqnarray*}
Then we readily obtain the asymptotic result for the test statistic $TCvM_n$ under $H_0$.

\begin{theorem}\label{Tcvm-null}
Suppose that Conditions (A1)-(A7) and (A9)-(A13) hold and $\Gamma^{tr}_{e}(z)$ is nonsingular for each $z \in \mathbb{R}$.
Under $H_0$, if $p^3(\nu_n\log{n})^6 = o(n)$, then
$$ TCvM_n \longrightarrow \int_0^1 B(t)^2 dt, \quad   in \ distribution, $$
where $B(t)$ is the standard Brownian motion.
\end{theorem}

According to Theorem \ref{Tcvm-null}, the test $TCvM_n$ based on the martingale-transformed process $\hat{T}_n \hat{U}_n(t)$ is asymptotically distribution free and then its asymptotic critical values can be tabulated.
%Furthermore, the test $TCvM_n$ can only have the sensitivity rate of order $n^{-1/4}$, rather than $n^{-1/2}$. This efficiency loss is an inherent characteristic of the test based the martingale transformed residual empirical process, no matter whether the dimension of parameters is divergent or fixed as the sample size tends to infinity.
For significance levels $0.01, 0.05$, and $0.1$, the critical values for $TCvM_n$ are about $2.787, 1.656$, and $1.196$, respectively. More critical values can be found in \cite{klein2003} (Table C.6, Page 478) or \cite{weller2009} (Table 1, Page 748).
For $t_0$, as suggested by \cite*{stute1998b}, we choose the $99\%$ quantile of $\hat{F}_{\hat e} $ in practice.
According to Condition (A11), we suggest choosing $a_n = n^{-2}$, $-c_{1n} = \zeta_{0.01} - \log^{1/2}{n}$, and $c_{2n}=\zeta_{0.99} + \log^{1/2}{n}$ in the simulation studies, where $\zeta_{0.01}$ and $\zeta_{0.99}$ are the $1\%$ and $99\%$ quantiles of $\hat{F}_{\hat e} $, respectively.

\section{Power analysis}
In this section, we investigate the asymptotic properties of the processes $\hat{U}_n(t)$ and $\hat{T}_n \hat{U}_n(t)$ under the alternative hypotheses in diverging-dimension settings. Consider the following sequence of alternative hypotheses converging to the null at the rate $r_n = n^{-\alpha}$:
\begin{equation}\label{4.1}
H_{1n}:  Y_n=m(X, \beta_0) + r_n S(X) + \varepsilon
\end{equation}
where $\alpha \in [0, 1/2]$, $S(\cdot)$ is a real-valued and non-constant function satisfying $E[S(X)]=0$, $\mathbb{P}(S(X)=0) < 1$, and $E|S(X)|^k = O(1)$.  Note that $\alpha=0$ corresponds to the global alternative and $\alpha >0$ to local alternatives.

First, we discuss the asymptotic properties of $\hat{U}_n(t)$ under $H_1$ in diverging-dimension scenarios. Lemma S6 and Remark S2 in the Supplementary Material show that under $H_{1n}$ with $r_n = n^{-\alpha}$ and $\alpha \in (0, 1/2]$,
\begin{eqnarray}\label{4.2}
\hat{U}_n(t)
= U_n^1(t) + f_{\varepsilon}(t) \sqrt{n}r_n \{M^{\top} \Sigma^{-1} M_S - E[g_0(X)S(X)]\} + R_n^1(t),
\end{eqnarray}
where $M=E[g_0(X)\dot{m}(X, \beta_0)]$, $\Sigma=E[\dot{m}(X, \beta_0)\dot{m}(X, \beta_0)^{\top}]$, $M_S=E[S(X) \dot{m}(X, \beta_0)]$, $U_n^1(t)$ is given in (\ref{2.6}), and the remainder $R_n^1(t)$ satisfies $\sup_{t}|R_n^1(t)| = o_p(\sqrt{n}r_n)$.
Consequently, under $H_{1n}$ with $r_n= n^{-1/2}$,
\begin{eqnarray}\label{4.3}
\hat{U}_n(t) = U_n^1(t) + f_{\varepsilon}(t) \{M^{\top} \Sigma^{-1} M_S - E[g_0(X)S(X)]\} + R_n^1(t),
\end{eqnarray}
where $R_n^1(t) = o_p(1)$ uniformly in $t$.
Combining (\ref{4.2}) and (\ref{4.3}) with Theorem \ref{expan-Un}, we derive the asymptotic properties of the weighted residual empirical process $\hat{U}_n(t)$ under various alternative hypotheses. Write $G_n(t)= E[g_0(X)F_{\varepsilon}(t+ m(X, \tilde{\beta}_0) - m(X))]$ and $ A_n = M^{\top} \Sigma^{-1} M_S - E[g_0(X)S(X)]$.

\begin{theorem}\label{Un-limitalter}
Suppose that Conditions (A1)-(A7) hold. \\
{\rm (i)} Under $H_1$, if $p^3\log^6{n} = o(n)$, then
$$ n^{-1/2} \hat{U}_n(t) \longrightarrow G(t),  \quad in \ probability, $$
where $G(t)$ is a pointwise limit of $G_n(t)$. \\
{\rm (ii)} Under $H_{1n}$ with $r_n = n^{-\alpha}$ and $\alpha \in (0, 1/2)$, if $p^3 \log^6{n} = o(r_n^{-2})$, then
$$ (n r_n^2)^{-1/2} \hat{U}_n(t) \longrightarrow A f_{\varepsilon}(t), \quad in \ probability, $$
where $A$ is the limit of $A_n$. \\
{\rm (iii)} Under $H_{1n}$ with $r_n = n^{-1/2}$, if $p^3\log^6{n} = o(n)$, then
$$ \hat{U}_n(t) \longrightarrow U_{\infty}^1(t) + A f_{\varepsilon}(t), \quad in \ distribution, $$
in the Skorohod Space $D[-\infty,\infty]$, where $U_{\infty}^1(t)$ is a Gaussian process given in Theorem \ref{hatU-limitnull}.
\end{theorem}

By Theorem \ref{Un-limitalter} and the Extended Continuous Mapping Theorem, we readily obtain the limiting distributions of the test $CvM_n$ based on the process $\hat{U}_n(t)$ under various alternative hypotheses.

\begin{theorem}\label{CVM-limitalter}
Suppose that Conditions (A1)-(A7) hold and $G(t)$ and $A f_{\varepsilon}(t)$ are nonzero functions. \\
{\rm (i)} Under $H_1$, if $p^3\log^6{n} = o(n)$, then $CvM_n \to \infty$ in probability at the rate of $n^{1/2}$. \\
{\rm (ii)} Under $H_{1n}$ with $r_n = n^{-\alpha}$ and $\alpha \in (0, 1/2)$, if $p^3 \log^6{n} = o(r_n^{-2})$, then $CvM_n \to \infty$ in probability at the rate of $\sqrt{n} r_n$. \\
{\rm (iii)} Under $H_{1n}$ with $r_n = n^{-1/2}$, if $p^3 \log^6{n} = o(n)$, then
\begin{eqnarray*}
CvM_n \longrightarrow \int_{\mathbb{R}} |\rho^{-1} [U_{\infty}^1(t) + Af_{\varepsilon}(t)]|^2 dF_{\varepsilon}(t),\quad in \ distribution.
\end{eqnarray*}
\end{theorem}

According to Theorems \ref{Un-limitalter} and \ref{CVM-limitalter}, if we choose a weight function $g(x)$ such that the functions $G(t)$ and $A f_{\varepsilon}(t)$ are nonzero, then the test $CvM_n$ based the process $\hat{U}_n(t)$ could be consistent with asymptotic power $1$ under the alternatives and can detect local alternatives converging to the null at a parametric rate $n^{-1/2}$.
Furthermore, the following result shows that the smoothing residual bootstrap is asymptotically consistent under various alternatives in diverging-dimension scenarios.

\begin{theorem}\label{bootstrap-alter}
Assume that Conditions (A1)-(A8) hold. \\
{\rm (i)} If $p^3 \log^6{n} = o(r_n^{-2})$, then, conditionally on the original sample $\mathcal{Y}_n$, under $H_{1n}$ with $r_n = n^{-\alpha}$ and $\alpha \in (0, 1/2]$, the results in Theorem \ref{bootstrap-null} continue to hold. \\
{\rm (ii)} If $p^3\log^6{n} = o(n)$, then, conditionally on the original sample $\mathcal{Y}_n$, under $H_1$, $CvM_n^*$ converges to a finite, weak limit which may be different from the weak limit of $CvM_n$ under $H_0$. This concludes that $ \lim_{n \to \infty} \mathbb{P}( CvM_n > c_{n, \zeta}^*) = 1 $.
\end{theorem}

However, interestingly and surprisingly, the test statistic $TCvM_n$ based on the martingale-transformed process $\hat{T}_n \hat{U}_n(t)$ may not achieve the parametric sensitive rate under the local alternatives.
According to (\ref{2.6}) and (\ref{4.3}), although the shift terms in the expansion of $\hat{U}_n(t)$ are different under $H_0$ and $H_{1n}$ with $r_n=n^{-1/2}$, they share the same shift functions $F_{\varepsilon}(t)$ and $f_{\varepsilon}(t)$ in both cases.
Note that one purpose of the martingale transformation $T_{\varepsilon}$ is to eliminate the functions $F_{\varepsilon}(t)$ and $f_{\varepsilon}(t)$ in the expansion of $\hat{U}_n(t)$.  Lemma~S8 in the Supplementary Material shows that under some regularity conditions and $H_{1n}$ with $r_n = n^{-1/2}$,
$$\hat{T}_nF_{\varepsilon}(t)= o_p(1), \ \ \hat{T}_n f_{\varepsilon}(t) = o_p(1), \ \ {\rm and} \ \
\hat{T}_n \hat{U}_n(t) = T_{\varepsilon}^{tr}U_n^0(t) + o_p(1), $$
uniformly in $t$.
Consequently, the martingale-transformed process $\hat{T}_n \hat{U}_n(t)$ has the same limiting distribution under both $H_0$ and $H_{1n}$ with $r_n=n^{-1/2}$.
In other words, the test statistic $TCvM_n$ based on $\hat{T}_n \hat{U}_n(t)$ {\it cannot} detect the local alternatives distinct from the null at the rate of order $n^{-1/2}$. Therefore, this test would lose power in theory compared to that based on the process $\hat{U}_n(t)$ without martingale transformation.
This is very different from existing asymptotically distribution-free tests based on martingale transformations of residual-marked empirical processes, which usually has the sensitivity rate of order $n^{-1/2}$, see \cite*{stute1998b}, \cite{bai2003}, and \cite{tan2019a} for instance.
We also note that when testing for membership in a location class of distributions, the asymptotically distribution-free test based on the martingale-transformed residual empirical process can still have the sensitivity rate of order $n^{-1/2}$ to local alternatives and even have higher power than its untransformed counterparts (\cite{khmaladze2009}).
This is because those existing martingale transformations in the literature usually do not eliminate the shift terms arising from the alternatives with $r_n= n^{-1/2}$.

To further check the exact rate of sensitivity of the transformed process $\hat{T}_n \hat{U}_n(t)$ to the local alternatives, we impose some additional conditions.

(A14). Under $H_{1n}$ with $r_n = n^{-\alpha}$ and $\alpha \in (0, 1/4]$, we have $\sup_{t \leq -c_{1n}} \|W_{\varepsilon}(t) \| = o_p(r_n)$ and $ r_n \sup_{ -c_{1n} \leq v \leq c_{2n}}|\frac{\dot{f}_{\varepsilon}(v)}{f_{\varepsilon}(v)}| = o(1)$, where $W_{\varepsilon}(t) = (F_{\varepsilon}(t), f_{\varepsilon}(t))^{\top}$.

The next theorem provides the exact sensitivity rate of the martingale-transformed process $\hat{T}_n \hat{U}_n(t)$ under the alternatives in diverging-dimension settings. The proof of this result will be provided in the Supplementary Material. Write $ T^{tr}_e G_n(t) = G_n(t) - \int_{-\infty}^t h^{tr}_e(z)^{\top} \Gamma^{tr}_e(z)^{-1} \int_z^{\infty} h^{tr}_e(v)dG_n(v) dF_e(z) $, where $G_n(t)= E[g_0(X)F_{\varepsilon}(t+ m(X, \tilde{\beta}_0) - m(X))]$.

\begin{theorem}\label{limit-TnUn-alter}
Suppose that Conditions (A1)-(A7) and (A9)-(A14) hold and $\Gamma^{tr}_{e}(z)$ is nonsingular for each $z \in \mathbb{R}$. \\
{\rm (i)} Under $H_1$, if $p^3(\nu_n\log{n})^6 = o(n)$, then
$$ n^{-1/2} \hat{T}_n \hat{U}_n(t) \longrightarrow L_1(t), \quad in \ probability, $$
where $L_1(t)$ is the pointwise limit of $T^{tr}_e G_n(t)$. \\
{\rm (ii)} Under $H_{1n}$ with $r_n = n^{-\alpha}$ and $\alpha \in (0, 1/4)$, if $p^3 (\nu_n \log{n})^6 = o(r_n^{-2})$, $\kappa_n^{-1} r_n^2 = o(1)$, and $\sqrt{p}r_n h^{-1} = o(1)$, then
$$ (n r_n^4)^{-1/2} \hat{T}_n \hat{U}_n(t) \longrightarrow L_2(t), \quad in \ probability, $$
where $L_2(t)$ is a determinist function that is specified in (18) in the Supplementary Material. \\
{\rm (iii)} Under $H_{1n}$ with $r_n = n^{-1/4}$, if $p^3 (\nu_n\log{n})^6 = o(\sqrt{n})$, then
$$ \hat{T}_n \hat{U}_n(t) \longrightarrow \rho B(F_{\varepsilon}(t)) + L_2(t), \quad in \ distribution, $$
in the Skorohod Space $D[-\infty,\infty)$.
\end{theorem}

It follows from Theorem \ref{limit-TnUn-alter} that under $H_{1n}$ with $r_n = n^{-1/4}$, the martingale-transformed process $\hat{T}_n \hat{U}_n(t)$ weakly converges to a Gaussian process that may be different from the limiting null distribution of $\hat{T}_n \hat{U}_n(t)$. Under $H_{1n}$ with $r_n = n^{-\alpha}$ and $\alpha \in (0, 1/4)$, the process $\hat{T}_n \hat{U}_n(t) $ may diverge to infinity in probability. Therefore, tests based on the martingale-transformed process $\hat{T}_n \hat{U}_n(t)$ may lose power dramatically, detecting local alternatives distinct from the null only at the rate of order $n^{-1/4}$.
Note that this efficiency loss is an inherent characteristic of tests based on the martingale-transformed residual empirical process, regardless of whether the dimension of parameters is fixed or divergent as the sample size $n$ tends to infinity. However, we also observe that this slower rate remains independent of the dimension of the parameter or predictor vector. Consequently, in high-dimensional scenarios, such tests remain useful and competitive compared to local smoothing tests  whose sensitivity rates to local alternatives are often dimension-dependent (e.g., \cite{hardle1993}, \cite{zheng1996}) as we commented in the Introduction section.

The following asymptotic result for the test statistic $TCvM_n$ is a consequence of Theorem \ref{limit-TnUn-alter} and the Extended Continuous Mapping Theorem.

\begin{theorem}\label{Tcvm-alter}
Suppose that Conditions (A1)-(A7) and (A9)-(A13) hold, $\Gamma^{tr}_{e}(z)$ is nonsingular for each $z \in \mathbb{R}$, and $L_1(t)$ and $L_2(t)$ are nonzero functions. \\
{\rm (i)} Under $H_1$, if $p^3(\nu_n\log{n})^6 = o(n)$, then $ TCvM_n \to \infty$ in probability at the rate of $\sqrt{n}$. \\
{\rm (ii)} Under $H_{1n}$ with $r_n= n^{-\alpha}$ and $\alpha \in (0, 1/4)$, if $p^3(\nu_n\log{n})^6 = o(r_n^{-2})$, $\kappa_n^{-1} r_n^2 = o(1)$, and $\sqrt{p}r_n h^{-1} = o(1)$, then $ TCvM_n \to \infty$ in probability at the rate of $\sqrt{n}r_n^2$. \\
{\rm (iii)} Under $H_{1n}$ with $r_n = n^{-1/4}$, if $p^3(\nu_n\log{n})^6 = o(\sqrt{n})$, then
$$ TCvM_n \longrightarrow \int_0^1 \big| B(t) + L_2(F_{\varepsilon}^{-1}(t F_{\varepsilon}(t_0))) / \sqrt{F_{\varepsilon}(t_0)} \big|^2 dt, \quad   in \ distribution. $$
\end{theorem}

\begin{remark}
The weighted residual empirical process and its martingale transformation can also be applied to assess whether the conditional variance function $\sigma^2(\cdot)$ in~(\ref{1.1}) belongs to a specific parametric family of functions $\tilde{\mathcal{M}}=\{ \sigma^2(\cdot, \theta): \theta \in \tilde{\Theta} \subset \mathbb{R}^q\}$, when the mean function is specified. The asymptotic properties of the weighted residual empirical process and its martingale transformation for testing conditional variance functions are very similar to those developed in Sections 2-4. For brevity, we omit these details in the main text and provide them in Sections S6-S7 of the Supplementary Material.
Additionally, it is worth noting that our methodology can be extended to heteroscedastic cases and dependent data. Further discussion on these extensions is deferred to Section S8 of the Supplementary Material to save space.
\end{remark}

\section{The choice of the weight function}
For practical use, we must delicately choose the weight function $g(\cdot)$ such that the tests based on $\hat{U}_n$ and $\hat{T}_n \hat{U}_n(t)$ can have good power.
Recall that under $H_0$, the assertion (\ref{2.6}) shows that
$$ \hat{U}_n(t) = U_n^1(t) + o_p(1), $$
uniformly in $t$.
Under $H_{1n}$ with $r_n= n^{-1/2}$, it follows from (\ref{4.3}) that
$$ \hat{U}_n(t) = U_n^1(t) + f_{\varepsilon}(t) \{M^{\top} \Sigma^{-1} M_S - E[g_0(X)S(X)]\} + o_p(1), $$
uniformly in $t$.
Therefore, to have good power performance, a natural idea is to choose a weight function $g(\cdot)$ such that  $ | M^{\top} \Sigma^{-1} M_S - E[g_0(X)S(X)] | $ is positive and is as large as possible.
Recall that $M=E[g_0(X) \dot{m}(X, \beta_0)]$ and $g_0(X) = g(X) - E[g(X)]$, it follows from the Cauchy-Schwarz inequality that
\begin{eqnarray*}
\rho_n^{-1} |M^{\top} \Sigma^{-1} M_S - E[g_0(X)S(X)]|
&=&    |E\{\frac{g_0(X)}{var\{g(X)\}^{1/2}} [\dot{m}(X, \beta_0)^{\top} \Sigma^{-1} M_S - S(X)] \}| \\
&\leq& var^{1/2}\{\dot{m}(X, \beta_0)^{\top} \Sigma^{-1} M_S - S(X) \}
\end{eqnarray*}
with the equality holds if $g(X)=\dot{m}(X, \beta_0)^{\top} \Sigma^{-1} M_S - S(X)$, where $\rho_n^2 = var\{g(X)\}$. Since $M_S=E[S(X) \dot{m}(X, \beta_0)]$ and $S(X) = r_n^{-1} \{m(X) - m(X, \beta_0)\}$ under $H_{1n}$, the theoretically optimal weight function under $H_{1n}$ should be
\begin{eqnarray}\label{6.1}
g^{opt}(X) = \dot{m}(X, \beta_0)^{\top} \Sigma^{-1} E\{\dot{m}(X, \beta_0) [m(X) - m(X, \beta_0)] \} - [m(X)- m(X, \beta_0)].
\end{eqnarray}
Here, the coefficient $r_n$ is dropped off as the weight function $g_0(X)$ is standardized in constructing the test statistic.

Although $g^{opt}(\cdot)$ is theoretically optimal under $H_{1n}$, it cannot be used as a weight function for practical use. Note that $m(X)= m(X, \beta_0)$ under $H_0$ and then $g^{opt}(X) \equiv 0$.
Consequently, the process $\hat{U}_n(t) \equiv 0$ for all $t \in \mathbb{R}$, which causes the theoretical results under the null to become meaningless.
Also note that the power of our tests increases along with the correlation between the weight function $g(X)$ and $\dot{m}(X, \beta_0)^{\top} \Sigma^{-1} M_S - S(X)$. This observation suggests choosing a weight function that is non-constant and highly correlated with the function $\dot{m}(X, \beta_0)^{\top} \Sigma^{-1} M_S - S(X)$. Therefore, we consider using $m(X)$ instead of $m(X) - m(X, \beta_0)$ in~(\ref{6.1}) and suggest a weight function as
\begin{equation}\label{6.2}
g(X) = \dot{m}(X, \beta_0)^{\top} \Sigma^{-1} E[\dot{m}(X, \beta_0) m(X)] - m(X).
\end{equation}
Furthermore, we check the performance of the resulting tests numerically based on the weight $g(X)$ and the optimal weight $g^{opt}(X)$ under the alternatives when the regression function $m(\cdot)$ is given.
These unreported simulations show that the test based on $g(X)$ is only slightly less powerful than that based on $g^{opt}(X)$. Therefore, we recommend the weight function $g(X)$ in~(\ref{6.2}) for practical use when the unknowns in $g(\cdot)$ are substituted with the corresponding approximations.

The regression function $m(X)$ in~(\ref{6.2}) is unknown under $H_1$ and has to be estimated in a nonparametric way. This introduces significant challenges in diverging-dimension settings due to the curse of dimensionality.
To address these challenges, we propose a novel procedure for estimating $m(x)$ based on Fourier transformation, which is inspired by \cite{stute2008}.
Let $L^2(F_X)$ be the Hilbert Space of squared integral functions endowed with the inner product $\left<g, h \right> = \int g(x)h(x)dF_X$ and let
$ \dot{m}_j(x, \tilde{\beta}_0) = \frac{\partial m(x, \tilde{\beta}_0)}{\partial \beta_j}$ for $ j=1, \cdots, p$, where $F_X$ denotes the cumulative distribution function of $X$.
Recall that $m(x)$ satisfies  $ E\{[m(X) - m(X, \tilde{\beta}_0)]\dot{m}(X, \tilde{\beta}_0)\}=0 $.
This means that $m(x) - m(x, \tilde{\beta}_0)$ is orthogonal to the span of $\{ \dot{m}_1, \cdots, \dot{m}_p \}$.
We then apply the Gram-Schmidt orthonormalization procedure to $\dot{m}_1, \cdots, \dot{m}_p$ and expand it to be an orthonormal basis $w_1, \cdots, w_p, w_{p+1}, w_{p+2}, \cdots$ of $L^2(F_X)$ with ${\rm span}\{w_1, \cdots, w_p \} = {\rm span} \{\dot{m}_1, \cdots, \dot{m}_p\}$. Consequently, $m(x)$ admits a Fourier representation
$$ m(x)= m(x, \tilde{\beta}_0) + \sum_{i=p+1}^{\infty} c_i w_i(x), $$
where $c_i = \langle m(x)-m(x, \tilde{\beta}_0), w_i \rangle$. To estimate the Fourier coefficients $c_i$, we observe
\begin{eqnarray*}
c_i = E\{[m(X)-m(X, \tilde{\beta}_0)]w_i(X)\} = E\{[Y-m(X, \tilde{\beta}_0)]w_i(X)\}.
\end{eqnarray*}
Their empirical analogues are $\hat{c}_i = n^{-1} \sum_{j=1}^n \hat{e}_j w_i(X_j)$ where $\hat{e}_j = Y_j-m(X_j, \hat{\beta})$. Then the estimator of $m(x)$ is given by $\hat{m}(x) = m(x, \hat{\beta}) + \sum_{i=p+1}^{\infty} \hat{c}_i w_i(x)$.
In practice, $\hat{m}(x)$ could be approximated by
\begin{eqnarray}\label{6.4}
\hat{m}_l(x) = m(x, \hat{\beta}) + \sum_{i=p+1}^{p+l} \hat{c}_i w_i(x),
\end{eqnarray}
for some fixed number $l$. Thus, our analysis for the nonparametric alternative leads to choosing a weight function as
\begin{equation}\label{6.5}
\hat{g}_l(X) = \dot{m}(X, \hat{\beta})^{\top} \hat{\Sigma}^{-1} \hat{E}[\dot{m}(X, \hat{\beta}) \hat{m}_l(X)] - \hat{m}_l(X),
\end{equation}
where $\hat{\Sigma}= n^{-1}\sum_{i=1}^n \dot{m}(X_i, \hat{\beta})\dot{m}(X_i, \hat{\beta})^{\top}$ and $\hat{E}[\dot{m}(X, \hat{\beta}) \hat{m}_l(X)] = n^{-1} \sum_{i=1}^n \dot{m}(X_i, \hat{\beta}) \hat{m}_l(X_i)$.
The number $l$ and the choice of orthonormal basis $w_{p+1}, w_{p+2}, \cdots, w_{p+l}$ (more precisely, the set of linearly independent functions) would depend on researchers.

Note that the estimated weight function $\hat{g}_l$ given in (\ref{6.5}) only involves parametric estimators. Following the same line as the proof for Theorem~{\ref{expan-Un}}, we can show that under $H_0$ and some regularity conditions,
\begin{eqnarray}\label{estiU}
\hat{U}_n(t) = \frac{1}{\sqrt{n}} \sum_{i=1}^n [\hat{g}_l(X_i) - \bar{\hat{g}}_l] I(\hat{e}_i \leq t) = \frac{1}{\sqrt{n}} \sum_{i=1}^n [g_l(X_i) - \bar{g}_l] I(\hat{e}_i \leq t) + R_n^0(t),
\end{eqnarray}
uniformly in $t$, where $\bar{\hat{g}}_l = n^{-1}\sum_{i=1}^n \hat{g}_l(X_i)$, $g_l(x) = \dot{m}(x, \beta_0)^{\top} \Sigma^{-1} E[\dot{m}(X, \beta_0) m_l(X)] +m_l(x) $ with $ m_l(x) = m(x, \beta_0) + \sum_{i=p+1}^{p+l} c_i w_i(x)$, and $\sup_t|R_n^0(t)| = o_p(1)$ under the condition $p^3\log^6{n} = o(n)$ and $\sup_t|R_n^0(t)| = o_p(\nu_n^{-1})$ under the condition $p^3(\nu_n\log{n})^6 = o(n)$.
The proof for~(\ref{estiU}) is provided in Lemma S10 of the Supplementary Material. Therefore, the result of Theorem \ref{hatU-limitnull} continues to hold when using the estimated weight function $\hat{g}_l$ in $\hat{U}_n(t)$, and the martingale transformation can also be applied to this process $\hat{U}_n(t)$. This implies that our tests based on $\hat{U}_n(t)$ with estimated weight function $\hat{g}_l$ and its martingale transformation $\hat{T}_n \hat{U}_n(t)$ can be applied in practice, although such a choice of the weight function $\hat{g}_l$ may not make the resulting tests omnibus.
Furthermore, we suggest using sufficient dimension reduction techniques to determine an appropriate orthonormal basis in this paper. More details will be provided in the next section.

\section{Numerical studies}
\subsection{Simulations}
In this subsection, we evaluate the performances of the tests $CvM_n$ and $TCvM_n$ in finite samples.
For the smoothing parameter $v_n$ in the smooth residual bootstrap for $CvM_n$, we choose $v_n= 0.2$ as suggested by \cite{dette2007}.
The sample $\{Z_i\}_{i=1}^n$ in the smoothing residual bootstrap are i.i.d. standard normal random variables.
The standard normal density function is also used as the kernel for the martingale-transformed process-based test $TCvM_n$. We compare our tests with \cite{zheng1996}'s test $T_n^{ZH}$, \cite{escanciano2006a}'s test $PCvM_n$, and \cite{xu2021}'s test $dCov_n$, although these tests are developed in fixed-dimension settings. The wild bootstrap is used to determine the critical values of $T_n^{ZH}$ and $PCvM_n$, while the residual bootstrap is used for $dCov_n$.
From the theoretical viewpoint in this paper, we consider $p=[3n^{1/3}] - 3$ with the sample sizes $n=100, 200, 400, 600$. To give a relatively thorough comparison among these tests, we also consider fixed-dimension cases with $p = 2,4,8$ and $n = 100$ in the simulation studies.  The significance level is set to be $\zeta = 0.05$. The simulation results are based on the averages of 1000 replications and the bootstrap approximation of $B = 500$ replications. In the following simulation studies, $a=0$ corresponds to the null while $a \neq 0$ to the alternatives.

{\it Study 1.} We generate data from the following regression models:
\begin{eqnarray*}
H_{11}:  Y_i &=& \beta_0^{\top}X_i + a (\beta_0^{\top}X_i)^2 + \varepsilon_i; \\
H_{12}:  Y_i &=& \beta_0^{\top}X_i + a \cos(0.6\pi \beta_0^{\top}X_i) + \varepsilon_i;
\end{eqnarray*}
where $\beta_0=(1, 1, \cdots, 1)^{\top}/\sqrt{p}$ and $X_i=(X_{i1}, \cdots, X_{ip})^{\top}$ is $N(0, \Sigma_1)$ or $N(0, \Sigma_2)$ independent of the standard Gaussian error term $\varepsilon$. Here $\Sigma_1=I_p$ and $\Sigma_2=(1/2^{|i-j|})_{p \times p}$.

First, we need to choose the weight function $g(X)$ to apply our tests.
%For parametric testing cases, $g(X)$ is determined by~(\ref{6.3}). For nonparametric testing cases,
In this paper, we use sufficient dimension reduction techniques to determine an orthonormal basis that can provide a good approximation of $m(X)$.
Let $\mathcal{S}_{Y|X}$ be the central subspace of $Y$ with respect to $X$. Here, the central subspace $\mathcal{S}_{Y|X}$ is defined as the intersection of all subspaces ${\rm span}(A)$ such that $Y \hDash X|A^{\top}X$, where $\hDash$ means the statistical independence and ${\rm span}(A)$ is the subspace spanned by the columns of the matrix $A$. Under mild conditions, such a subspace $\mathcal{S}_{Y|X}$ always exists (e.g., \cite{cook2009}). If $\mathcal{S}_{Y|X} = {\rm span}(B)$, then we have $m(X) = E(Y|X)=E(Y|B^{\top}X)$, where $B = (B_1, \cdots, B_s)$ is a $p \times s $ orthonormal matrix with $s \leq p$ and $s$ is the structural dimension. If $s=p$, there is no dimension reduction structure in regression models. In this paper, we adopt cumulative slicing estimation (CSE, \cite{zhu2010b}) to identify the central subspace $\mathcal{S}_{Y|X}$ and the minimum ridge-type eigenvalue ratio (MRER, \cite{zhux2017}) estimator to identify the structural dimension $s$, as they allow the divergence rate of $p$ to be  $p = o(\sqrt{n})$ and are very easy to implement in practice.
Let $\hat{B}=(\hat{B}_1, \cdots, \hat{B}_{\hat{s}})$ be a consistent estimator of $B$ obtained by CSE and MRER. Note that $\dot{m}(x, \tilde{\beta_0}) = (x_1, \cdots, x_p)^{\top}$ in Study 1. The functions selected for the Gram-Schimdt procedure to obtain the orthonormal basis are $\{ (\hat{B}_i^{\top}x)^2, (\hat{B}_i^{\top}x)^3, (\hat{B}_i^{\top}x)^4,  1 \leq i \leq \hat{s} \}$. If $\hat{m}(x)$ is determined by this orthonormal basis, then the weight function $g(x)$ is given by~(\ref{6.5}).
Another issue is the choice of the bandwidth $h$ for the test $TCvM_n$ based on the martingale-transformed process $\hat{T}_n\hat{U}_n(t)$. In this paper, we adopt cross-validation to select the bandwidth automatically.

The simulation results of Study 1 are presented in Tables 1-2. It can be observed that $CvM_n, TCvM_n$ and $PCvM_n$ can control the empirical size very well across all models and dimension settings. While $T_n^{ZH}$ and $dCov_n$ are usually conservative with smaller empirical sizes.
The test $TCvM_n$ based on the martingale-transformed process usually loses powers compared with the corresponding test without martingale transformation. This phenomenon aligns with the theoretical results in Section 4. For the low-frequency model $H_{11}$, the empirical powers of $CvM_n$ and $PCvM_n$ increase rapidly. However, for the high frequency model $H_{12}$, $PCvM_n$ and $TCvM_n$ have much lower power than $CvM_n$.
The tests $T_n^{ZH}$ and $dCov_n$ deteriorate rapidly in all models and dimension cases as the dimension $p$ increases.
Interestingly, the tests $CvM_n, TCvM_n$, and $PCvM_n$ appear to be less affected by the dimensionality.

\begin{table}[ht!]\caption{Empirical sizes and powers of the tests for $H_{11}$ in Study 1.}
\centering
{\small\scriptsize\hspace{8cm}
\renewcommand{\arraystretch}{0.6}\tabcolsep 0.4cm
\begin{tabular}{*{20}{c}}
\hline
&\multicolumn{1}{c}{a}&\multicolumn{1}{c}{n=100}&\multicolumn{1}{c}{n=100}&\multicolumn{1}{c}{n=100}&\multicolumn{1}{c}{n=100}&\multicolumn{1}{c}{n=200}&\multicolumn{1}{c}{n=400}&\multicolumn{1}{c}{n=600}\\
&&\multicolumn{1}{c}{p=2}&\multicolumn{1}{c}{p=4}&\multicolumn{1}{c}{p=8}&\multicolumn{1}{c}{p=10}&\multicolumn{1}{c}{p=14}&\multicolumn{1}{c}{p=19}&\multicolumn{1}{c}{p=22}\\
\hline
$CvM_n, \ \Sigma_1$      &0.00 &0.048   &0.051   &0.050  &0.065 &0.056   &0.060   &0.068\\
                         &0.15 &0.338   &0.332   &0.266  &0.264 &0.576   &0.905   &0.986\\
%\hline
$TCvM_n, \ \Sigma_1$     &0.00 &0.058   &0.064   &0.037   &0.047
                               &0.059   &0.055   &0.056\\
                         &0.15 &0.201   &0.160   &0.149   &0.162
                               &0.196   &0.257   &0.308\\
%\hline
$PCvM_n, \ \Sigma_1$     &0.00 &0.054   &0.068   &0.059  &0.074 &0.056   &0.047   &0.053\\
                         &0.15 &0.267   &0.290   &0.312  &0.298 &0.506   &0.805   &0.927\\
%\hline
$T_n^{ZH}, \ \Sigma_1$   &0.00 &0.045   &0.059   &0.048  &0.043 &0.038   &0.026   &0.029\\
                         &0.15 &0.111   &0.086   &0.052  &0.050 &0.047   &0.033   &0.050\\

$dCov_n, \ \Sigma_1$     &0.00 &0.052   &0.042   &0.030  &0.036 &0.032   &0.030   &0.037\\
                         &0.15 &0.157   &0.090   &0.057  &0.041 &0.044   &0.063   &0.100\\
\hline
\\
\hline
$CvM_n, \ \Sigma_2$      &0.00 &0.045   &0.053   &0.035  &0.054  &0.071   &0.057   &0.055\\
                         &0.15 &0.651   &0.865   &0.934  &0.955  &1.000   &1.000   &1.000\\
%\hline
$TCvM_n, \ \Sigma_2$     &0.00 &0.048   &0.059   &0.064   &0.049
                               &0.054   &0.060   &0.062\\
                         &0.15 &0.304   &0.405   &0.408   &0.410
                               &0.587   &0.870   &0.959\\
%\hline
$PCvM_n, \ \Sigma_2$     &0.00 &0.070   &0.051   &0.058  &0.069  &0.054   &0.061   &0.048\\
                         &0.15 &0.504   &0.768   &0.875  &0.917  &0.996   &1.000   &1.000\\
%\hline
$T_n^{ZH}, \ \Sigma_2$   &0.00 &0.060   &0.041   &0.047  &0.043  &0.047   &0.032   &0.029\\
                         &0.15 &0.239   &0.238   &0.145  &0.124  &0.086   &0.080   &0.091\\

$dCov_n, \ \Sigma_2$     &0.00 &0.042   &0.045   &0.044  &0.042  &0.046   &0.039   &0.043\\
                         &0.15 &0.466   &0.605   &0.541  &0.497  &0.826   &0.993   &1.000\\
\hline
\end{tabular}}
\end{table}

\begin{table}[ht!]\caption{Empirical sizes and powers of the tests for $H_{12}$ in Study 1.}
\centering
{\small\scriptsize\hspace{8cm}
\renewcommand{\arraystretch}{0.6}\tabcolsep 0.4cm
\begin{tabular}{*{20}{c}}
\hline
&\multicolumn{1}{c}{a}&\multicolumn{1}{c}{n=100}&\multicolumn{1}{c}{n=100}&\multicolumn{1}{c}{n=100}&\multicolumn{1}{c}{n=100}&\multicolumn{1}{c}{n=200}&\multicolumn{1}{c}{n=400}&\multicolumn{1}{c}{n=600}\\
&&\multicolumn{1}{c}{p=2}&\multicolumn{1}{c}{p=4}&\multicolumn{1}{c}{p=8}&\multicolumn{1}{c}{p=10}&\multicolumn{1}{c}{p=14}&\multicolumn{1}{c}{p=19}&\multicolumn{1}{c}{p=22}\\
\hline
$CvM_n, \ \Sigma_1$      &0.00  &0.048   &0.052   &0.059  &0.056   &0.056   &0.063   &0.059\\
                         &0.75  &0.896   &0.842   &0.693  &0.623   &0.935   &1.000   &1.000\\
%\hline
$TCvM_n, \ \Sigma_1$     &0.00  &0.063   &0.056   &0.052   &0.047
                                &0.049   &0.051   &0.053\\
                         &0.75  &0.169   &0.155   &0.115   &0.095
                                &0.139   &0.269   &0.385\\
%\hline
$PCvM_n, \ \Sigma_1$     &0.00  &0.066   &0.059   &0.069  &0.062   &0.050   &0.057   &0.045\\
                         &0.75  &0.306   &0.253   &0.257  &0.237   &0.380   &0.619   &0.793\\
%\hlin
$T_n^{ZH}, \ \Sigma_1$   &0.00  &0.041   &0.043   &0.047  &0.041   &0.048   &0.038   &0.025\\
                         &0.75  &0.927   &0.616   &0.160  &0.101   &0.094   &0.055   &0.046\\

$dCov_n, \ \Sigma_1$     &0.00  &0.033   &0.056   &0.030  &0.033   &0.022   &0.032   &0.029\\
                         &0.75  &0.782   &0.318   &0.071  &0.062   &0.076   &0.099   &0.111\\
\hline
\\
\hline
$CvM_n, \ \Sigma_2$      &0.00  &0.046   &0.047   &0.048 &0.057   &0.054   &0.055   &0.067\\
                         &0.75  &0.870   &0.670   &0.396 &0.335   &0.436   &0.507   &0.585\\
%\hline
$TCvM_n, \ \Sigma_2$     &0.00  &0.051   &0.043   &0.062   &0.056
                                &0.055   &0.059   &0.069\\
                         &0.75  &0.170   &0.130   &0.087   &0.090
                                &0.065   &0.060   &0.057\\
%\hline
$PCvM_n, \ \Sigma_2$     &0.00  &0.039   &0.059   &0.064 &0.064   &0.070   &0.053   &0.057\\
                         &0.75  &0.148   &0.067   &0.064 &0.063   &0.057   &0.067   &0.062\\
%\hline
$T_n^{ZH}, \ \Sigma_2$   &0.00  &0.052   &0.057   &0.044 &0.053   &0.051   &0.047   &0.034\\
                         &0.75  &0.936   &0.691   &0.210 &0.147   &0.110   &0.085   &0.077\\

$dCov_n, \ \Sigma_2$     &0.00  &0.046   &0.029   &0.035 &0.031   &0.040   &0.027   &0.038\\
                         &0.75  &0.823   &0.308   &0.071 &0.043   &0.043   &0.048   &0.060\\
\hline
\end{tabular}}
\end{table}

The hypothetical models in Study 1 are single-index models under both the null and alternatives. Next, we consider more complex models in the second simulation study.

{\it Study 2.} Generate data from the following models:
\begin{eqnarray*}
H_{21}:  Y_i &=& \beta_0^{\top}X_i + a \exp(\beta_1^{\top}X_i) + \varepsilon_i;  \\
H_{22}:  Y_i &=& \beta_0^{\top}X_i + a \{(\beta_1^{\top}X_i)^3 + \sin(0.5 \pi \beta_1^{\top}X_i) +
         (\beta_0^{\top}X_i) (\beta_1^{\top}X_i)\} + \varepsilon_i; \\
H_{23}:  Y_i &=& X_{i1} + a\{ |X_{i2}| + X_{i3}^3 - X_{i4}^2 + X_{i5}^3 + X_{i6}X_{i7} + \cos(\pi X_{i8}) + \sin(\pi X_{i9} X_{i10}) \} +
         \varepsilon_i,
\end{eqnarray*}
where $\beta_0=(1, 1, \dots, 1)^{\top}/\sqrt{p}$, $\beta_1=(0,\dots,0,\underbrace{1,\dots,1}_{p_1})^{\top} /\sqrt{p_1}$ with $p_1=[p/2]$,  and the covariate $X_i$ is the same as in Study 1. We do not consider the cases with $p=2, 4, 8$ and $n=100$ in model $H_{23}$, as the dimension of covariates is greater than or equal to $10$. Additionally, $H_{23}$ is a regression model without dimension reduction structure when $p=10$.

The choices for the weight function $g(x)$ and the bandwidth $h$ are the same as in Study 1.
We can observe that all tests can maintain the nominal level in most cases. For the model $H_{22}$, the empirical powers of $CvM_n, TCvM_n$, and $PCvM_n$ grow very quickly, and $CvM_n$ has better power performances than the other two. However, the test $PCvM_n$ performs much better than $CvM_n$ and $TCvM_n$ for model $H_{21}$. This may be due to the inaccurate estimation of the unknown regression function under the alternatives for our tests when the coefficient $a$ is too small.
For a more complicated model $H_{23}$, our tests $CvM_n, TCvM_n$, and $dCov_n$ have much higher empirical powers than $PCvM_n$ and $T_n^{ZH}$ in all cases. It is worth noting that $H_{23}$ is a regression model with no dimension reduction structure when $p=10$. This may suggest that, even though sufficient dimension reduction methods are used to choose the weight function, they can still be applied to check the adequacy of regression models without dimension reduction structure. Another interesting result is that $dCov_n$ generally exhibits much higher empirical powers in cases with correlated components of covariates compared to uncorrelated cases.

\begin{table}[ht!]\caption{Empirical sizes and powers of the tests for $H_{21}$ in Study 2.}
\centering
{\small\scriptsize\hspace{8cm}
\renewcommand{\arraystretch}{0.6}\tabcolsep 0.4cm
\begin{tabular}{*{20}{c}}
\hline
&\multicolumn{1}{c}{a}&\multicolumn{1}{c}{n=100}&\multicolumn{1}{c}{n=100}&\multicolumn{1}{c}{n=100}&\multicolumn{1}{c}{n=100}&\multicolumn{1}{c}{n=200}&\multicolumn{1}{c}{n=400}&\multicolumn{1}{c}{n=600}\\
&&\multicolumn{1}{c}{p=2}&\multicolumn{1}{c}{p=4}&\multicolumn{1}{c}{p=8}&\multicolumn{1}{c}{p=10}&\multicolumn{1}{c}{p=14}&\multicolumn{1}{c}{p=19}&\multicolumn{1}{c}{p=22}\\
\hline
$CvM_n, \ \Sigma_1$      &0.00 &0.055   &0.042   &0.041  &0.056   &0.055   &0.060   &0.052\\
                         &0.15 &0.143   &0.135   &0.122  &0.111   &0.188   &0.294   &0.418\\
%\hline
$TCvM_n, \ \Sigma_1$     &0.00 &0.048   &0.053   &0.049   &0.060
                               &0.057   &0.058   &0.058\\
                         &0.15 &0.134   &0.151   &0.119   &0.115
                               &0.166   &0.154   &0.189\\
%\hline
$PCvM_n, \ \Sigma_1$     &0.00 &0.047   &0.063   &0.053  &0.061   &0.059   &0.051   &0.045\\
                         &0.15 &0.655   &0.663   &0.667  &0.624   &0.924   &0.999   &0.999\\
%\hline
$T_n^{ZH}, \ \Sigma_1$   &0.00 &0.056   &0.050   &0.040  &0.032   &0.029   &0.023   &0.023\\
                         &0.15 &0.269   &0.153   &0.085  &0.067   &0.058   &0.040   &0.036\\

$dCov_n, \ \Sigma_1$     &0.10 &0.060   &0.036   &0.035  &0.024   &0.021   &0.037   &0.028\\
                         &0.15 &0.118   &0.103   &0.070  &0.071   &0.089   &0.101   &0.128\\
\hline
\\
\hline
$CvM_n, \ \Sigma_2$      &0.00 &0.050   &0.063   &0.047  &0.055   &0.054   &0.061   &0.061\\
                         &0.15 &0.167   &0.370   &0.537  &0.560   &0.896   &0.996   &1.000\\
%\hline
$TCvM_n, \ \Sigma_2$     &0.00 &0.056   &0.053   &0.061   &0.069
                               &0.051   &0.069   &0.052\\
                         &0.15 &0.169   &0.304   &0.377   &0.368
                               &0.526   &0.668   &0.874\\
%\hline
$PCvM_n, \ \Sigma_2$     &0.00 &0.052   &0.055   &0.059  &0.059   &0.055   &0.071   &0.053\\
                         &0.15 &0.630   &0.813   &0.917  &0.914   &0.998   &0.998   &1.000\\
%\hline
$T_n^{ZH}, \ \Sigma_2$   &0.00 &0.047   &0.056   &0.057  &0.059   &0.035   &0.031   &0.045\\
                         &0.15 &0.275   &0.232   &0.169  &0.138   &0.130   &0.155   &0.130\\

$dCov_n, \ \Sigma_2$     &0.00 &0.038   &0.044   &0.041  &0.034   &0.036   &0.051   &0.047\\
                         &0.15 &0.132   &0.307   &0.513  &0.566   &0.876   &0.994   &0.997\\
\hline
\end{tabular}}
\end{table}

\begin{table}[ht!]\caption{Empirical sizes and powers of the tests for $H_{22}$ in Study 2.}
\centering
{\small\scriptsize\hspace{8cm}
\renewcommand{\arraystretch}{0.6}\tabcolsep 0.4cm
\begin{tabular}{*{20}{c}}
\hline
&\multicolumn{1}{c}{a}&\multicolumn{1}{c}{n=100}&\multicolumn{1}{c}{n=100}&\multicolumn{1}{c}{n=100}&\multicolumn{1}{c}{n=100}&\multicolumn{1}{c}{n=200}&\multicolumn{1}{c}{n=400}&\multicolumn{1}{c}{n=600}\\
&&\multicolumn{1}{c}{p=2}&\multicolumn{1}{c}{p=4}&\multicolumn{1}{c}{p=8}&\multicolumn{1}{c}{p=10}&\multicolumn{1}{c}{p=14}&\multicolumn{1}{c}{p=19}&\multicolumn{1}{c}{p=22}\\
\hline
$CvM_n, \ \Sigma_1$      &0.00 &0.064   &0.052   &0.060  &0.057   &0.061   &0.056   &0.067\\
                         &0.25 &0.657   &0.625   &0.594  &0.541   &0.865   &1.000   &1.000\\
%\hline
$TCvM_n, \ \Sigma_1$     &0.00 &0.051   &0.052   &0.053   &0.054
                               &0.053   &0.056   &0.066\\
                         &0.25 &0.425   &0.400   &0.386   &0.329
                               &0.438   &0.575   &0.666\\
%\hline
$PCvM_n, \ \Sigma_1$     &0.00 &0.046   &0.047   &0.056  &0.056   &0.065   &0.066   &0.061\\
                         &0.25 &0.353   &0.341   &0.334  &0.328   &0.568   &0.855   &0.955\\
%\hline
$T_n^{ZH}, \ \Sigma_1$   &0.00 &0.059   &0.047   &0.050  &0.037   &0.040   &0.030   &0.025\\
                         &0.25 &0.426   &0.208   &0.081  &0.067   &0.068   &0.057   &0.037\\

$dCov_n, \ \Sigma_1$     &0.00 &0.038   &0.036   &0.031  &0.033   &0.032   &0.035   &0.032\\
                         &0.25 &0.541   &0.313   &0.172  &0.128   &0.206   &0.351   &0.497\\
\hline
\\
\hline
$CvM_n, \ \Sigma_2$      &0.00 &0.047   &0.042   &0.053  &0.048   &0.069   &0.065   &0.063\\
                         &0.25 &0.857   &0.968   &0.983  &0.994   &1.000   &1.000   &1.000\\
%\hline
$TCvM_n, \ \Sigma_2$     &0.00 &0.043   &0.063   &0.050   &0.065
                               &0.060   &0.051   &0.064\\
                         &0.25 &0.487   &0.624   &0.666   &0.697
                               &0.878   &0.953   &0.990\\
%\hline
$PCvM_n, \ \Sigma_2$     &0.00 &0.057   &0.067   &0.067  &0.068   &0.067   &0.056   &0.056\\
                         &0.25 &0.560   &0.711   &0.657  &0.610   &0.819   &0.960   &0.992\\
%\hline
$T_n^{ZH}, \ \Sigma_2$   &0.00 &0.057   &0.073   &0.049  &0.036   &0.048   &0.037   &0.042\\
                         &0.25 &0.568   &0.737   &0.611  &0.461   &0.561   &0.504   &0.498\\

$dCov_n, \ \Sigma_2$     &0.00 &0.039   &0.040   &0.028  &0.039   &0.033   &0.030   &0.027\\
                         &0.25 &0.683   &0.874   &0.929  &0.941   &0.999   &1.000   &1.000\\
\hline
\end{tabular}}
\end{table}

\begin{table}[ht!]\caption{Empirical sizes and powers of the tests for $H_{23}$ in Study 2.}
\centering
{\small\scriptsize\hspace{8cm}
\renewcommand{\arraystretch}{0.6}\tabcolsep 0.4cm
\begin{tabular}{*{20}{c}}
\hline
&\multicolumn{1}{c}{a}&\multicolumn{1}{c}{n=100}&\multicolumn{1}{c}{n=200}&\multicolumn{1}{c}{n=400}&\multicolumn{1}{c}{n=600}\\
&&\multicolumn{1}{c}{p=10}&\multicolumn{1}{c}{p=14}&\multicolumn{1}{c}{p=19}&\multicolumn{1}{c}{p=22}\\
\hline
$CvM_n, \ \Sigma_1$      &0.0   &0.054   &0.049   &0.065   &0.058\\
                         &1.0   &0.467   &0.795   &0.979   &0.997\\
%\hline
$TCvM_n, \ \Sigma_1$     &0.0   &0.058   &0.065   &0.056   &0.049\\
                         &1.0   &0.425   &0.560   &0.684   &0.767\\
%\hline
$PCvM_n, \ \Sigma_1$     &0.0   &0.055   &0.067   &0.054   &0.067\\
                         &1.0   &0.098   &0.111   &0.173   &0.261\\
%\hline
$T_n^{ZH}, \ \Sigma_1$   &0.0   &0.055   &0.038   &0.028   &0.024\\
                         &1.0   &0.303   &0.227   &0.148   &0.124\\

$dCov_n, \ \Sigma_1$     &0.0   &0.027   &0.033   &0.039   &0.033\\
                         &1.0   &0.658   &0.894   &0.993   &1.000\\
\hline
\\
\hline
$CvM_n, \ \Sigma_2$      &0.0   &0.051   &0.041   &0.048   &0.050\\
                         &1.0   &0.754   &0.958   &0.999   &1.000\\
%\hline
$TCvM_n, \ \Sigma_2$     &0.0   &0.058   &0.058   &0.060   &0.049\\
                         &1.0   &0.551   &0.697   &0.823   &0.902\\
%\hline
$PCvM_n, \ \Sigma_2$     &0.0   &0.060   &0.060   &0.051   &0.051\\
                         &1.0   &0.212   &0.308   &0.485   &0.593\\
%\hline
$T_n^{ZH}, \ \Sigma_2$   &0.0   &0.055   &0.036   &0.036   &0.039\\
                         &1.0   &0.609   &0.536   &0.427   &0.402\\

$dCov_n, \ \Sigma_2$     &0.0   &0.039   &0.042   &0.041   &0.046\\
                         &1.0   &0.869   &0.984   &1.000   &1.000\\
\hline
\end{tabular}}
\end{table}

It is noteworthy to describe the computational complexity when applying the test $CvM_n$ based the process $\hat{U}_n(t)$ and the test $TCvM_n$ based on the martingale-transformed process $\hat{T}_n \hat{U}_n(t)$. Here, we present the computational times from Study $H_{11}$ for illustration. These records are collected by running the tests in R with a Dell Precision 3530 Mobile Workstation with parallel computation. Table 6 records the time costs for 1000 replications in Study $H_{11}$ and 500 bootstrap sample. It is evident that the bootstrap-based test $CvM_n$ costs much more running time, while the test $TCvM_n$ based on the transformed process $\hat{T}_n \hat{U}_n(t)$ provides computational simplicity at the cost of power loss.

\begin{table}[ht!]\caption{The running time (secs) of 1000 replications for $H_{11}$ with $X \sim N(0, I_p)$.}
\centering
{\small\scriptsize\hspace{8cm}
\renewcommand{\arraystretch}{0.6}\tabcolsep 0.4cm
\begin{tabular}{*{20}{c}}
\hline
&\multicolumn{1}{c}{Time}&\multicolumn{1}{c}{n=100}&\multicolumn{1}{c}{n=100}&\multicolumn{1}{c}{n=100}&\multicolumn{1}{c}{n=100}&\multicolumn{1}{c}{n=200}&\multicolumn{1}{c}{n=400}&\multicolumn{1}{c}{n=600}\\
&&\multicolumn{1}{c}{p=2}&\multicolumn{1}{c}{p=4}&\multicolumn{1}{c}{p=8}&\multicolumn{1}{c}{p=10}&\multicolumn{1}{c}{p=14}&\multicolumn{1}{c}{p=19}&\multicolumn{1}{c}{p=22}\\
\hline
%&$CvM_n^d$        &38.48   &38.07   &40.17  &42.30   &82.86   &411.29   &892.06\\
%\hline
%&$TCvM_n^d$       &6.890   &7.300   &7.770  &8.040   &20.64   &84.10    &212.67\\
%\hline
&$CvM_n$          &48.48   &45.37   &46.16  &50.53   &106.55  &470.95   &1010.83\\
%\hline
&$TCvM_n$         &10.32   &10.57   &11.50  &12.27   &28.25   &104.50   &246.38\\
\hline
\end{tabular}}
\end{table}

\subsection{A real data example}
We apply the proposed tests to the Baseball Salaries Data Set which can be obtained from the website \url{https://www4.stat.ncsu.edu/~boos/var.select/baseball.html}. There are 337 Major League Baseball players salary information in this data set, with the salary $Y$ from the year 1992 as the output and $16$ performance measures from the year 1991 as predictor variables. The performance measures are $X_1$: Batting average, $X_2$: On-base percentage, $X_3$: the number of runs, $X_4$: the number of hits, $X_5$: the number of doubles, $X_6$: the number of triples, $X_7$: the number of home runs, $X_8$: the number of runs batted in, $X_9$: the number of walks, $X_{10}$: the number of strike-outs, $X_{11}$: the number of stolen bases, $X_{12}$: the number of errors, $X_{13}$: Indicator of free agency eligibility, $X_{14}$: Indicators of a free agent in 1991/2, $X_{15}$: Indicators of arbitration eligibility, and $X_{16}$: Indicators of arbitration in 1991/2. For easy interpretation, we standardize all variables respectively. To find out the regression relationship between the salary $Y$ and the performance measures $X=(X_1, \cdots, X_{16})^{\top}$, we first apply our tests to see whether a linear regression model $Y=\beta^{\top}X + \varepsilon$ is adequate to fit this data set. The values of test statistics are $CvM_n= 5.0851$ and $TCvM_n = 2.7918$, with the related $p$-values being about $0$ and $0.01$, respectively. This suggests rejecting the linearity of the underlying model. We further present a scatter plot of $Y$ against $\hat{\beta}_1^{\top}X$ in Figure 1(a) and the residuals against $\hat{\beta}_1^{\top}X$ in Figure 1(b), where $\hat{\beta}_1$ is the ordinary least squares estimator from the linear model. This plot also shows that a linear regression model for $(Y, X)$ may not be reasonable. From Figure 1(a), it seems that there may exist a quadratic relationship between $Y$ and $\hat{\beta}_1^{\top}X$. Therefore, we consider the following quadratic regression model to fit this data set
\begin{equation}\label{quadratic}
Y=a + \beta^{\top}X + c (\beta^{\top}X)^2 + \varepsilon.
\end{equation}
When applying our tests to model (\ref{quadratic}), the values of the test statistics are $CvM_n= 0.4298$ with a $p$-value $0.18$  and $TCvM_n = 0.5278$ with a $p$-value $0.31$, respectively. This suggests that the quadratic model (\ref{quadratic}) may be plausible. To further visualize this fit, Figure 2(a) presents a scatter plot of $Y$ against $\hat{\beta}_2^{\top}X$ with the fitted quadratic curve on the scatter plot, where $\hat{\beta}_2$ is the least squares estimator of $\beta$ obtained from model (\ref{quadratic}). Figure 2(b) is the scatter plot of residuals from the model (\ref{quadratic}) against the fitted values $\hat{Y}$, where $\hat{Y}= \hat{a} + \hat{\beta}_2^{\top}X + \hat{c} (\hat{\beta}_2^{\top}X)^2$. It is readily seen that there exists no trend between the residuals and the fitted values $\hat{Y}$. Therefore, the quadratic regression model (\ref{quadratic}) is adequate for this data set.

\begin{figure}
  \centering
  \includegraphics[width=13cm,height=5cm]{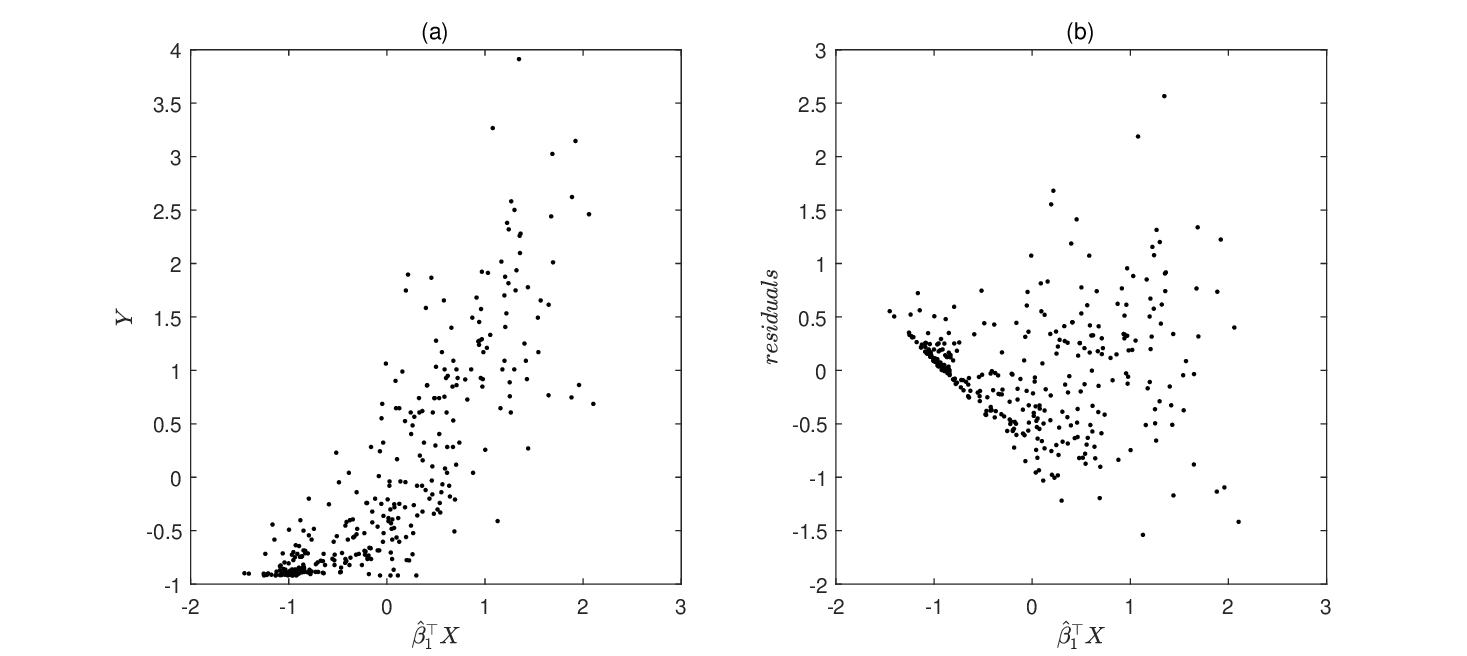}
  \caption{(a) Scatter plot of $Y$ versus $\hat{\beta}_1^{\top}X$ with $\hat{\beta}_1$ obtained from a linear model. (b) Scatter plot of $residuals$ from a linear model versus the fitted value $\hat{\beta}_1^{\top}X$. }\label{Figure1}
\end{figure}

\begin{figure}
  \centering
  \includegraphics[width=13cm,height=5cm]{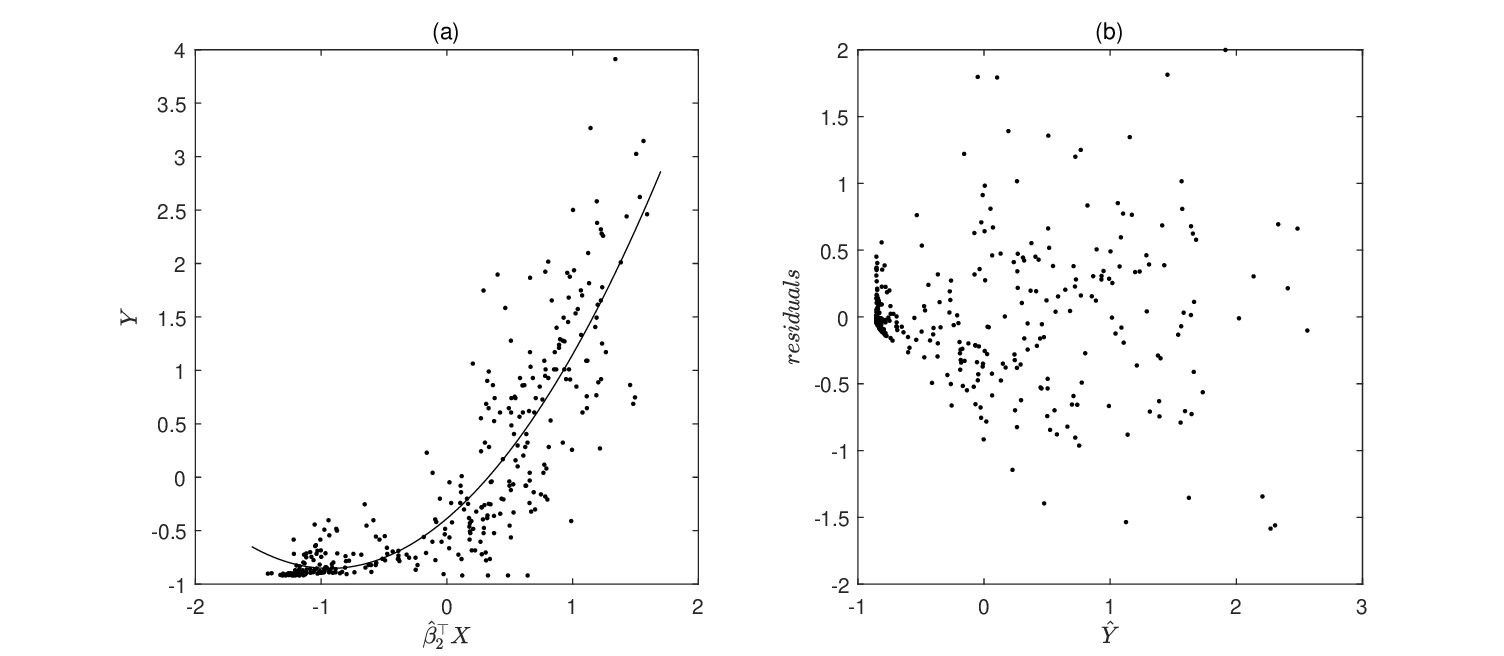}
  \caption{(a) Scatter plot of $Y$ versus $\hat{\beta}_2^{\top}X$ with $\hat{\beta}_2$ obtained from model (\ref{quadratic}). (b) Scatter plot of $residuals$ from the quadratic model (\ref{quadratic}) versus the fitted value $\hat{Y}$. }\label{Figure2}
\end{figure}

\section{Discussion}
This paper investigates the asymptotic properties of the weighted residual empirical process and its martingale transformation in diverging-dimension scenarios. These processes are used to construct goodness-of-fit tests for the mean and variance functions in regression models. The martingale transformation-based test is asymptotically distribution-free.
The test based on the weighted residual empirical process can detect local alternatives distinct from the null at the parametric rate of order $n^{-1/2}$. However, the test based on the martingale-transformed weighted residual empirical process loses power and can detect local alternatives distinct from the null only at the rate of order $n^{-1/4}$. This finding is interesting and somewhat surprising because existing martingale transformation-based tests in the literature typically can detect local alternatives departing from the null at the rate of order $n^{-1/2}$. This also implies, without considering computational issues, the test based on the weighted residual empirical process may dominate its martingale-transformed counterpart in power performance. Theoretical results are further validated through simulation studies.
As the test based on the weighted residual empirical process is not asymptotically distribution-free, a smooth residual bootstrap is used to approximate the limiting null distribution of this test. The validity of the bootstrap approximation is established in diverging-dimension settings. It is worth mentioning that, recently, \cite{neumeyer2019} proved that the classic residual bootstrap without smoothness is asymptotically valid for residual empirical processes in fixed-dimension cases. However, extending their method to our case with a diverging number of parameters may have some essential difficulties. This is because the empirical process $E_n$ involved in Lemma 4 of \cite{neumeyer2019} admits a diverging bracket number, making it challenging to obtain the asymptotic equicontinuity in diverging-dimension settings considered in the current paper. We also conducted simulations based on the classic residual bootstrap, which performs very well with large dimensions. Therefore, we conjecture that the classic residual bootstrap is still valid when the dimension diverges. %This is beyond the scope of this paper and deserves further study.

Note that model specification tests for the conditional mean and variance functions are special cases of testing conditional moment restrictions. It is of interest to  extend our method to test general (nonlinear in variables) conditional moment restrictions in high dimensional settings. This is beyond the scope of the current paper and deserves further study. Furthermore, as with many works (e.g., \cite{keilegom2008}, \cite{stute2008}, \cite{escanciano2010}, \cite{chown2018}), our methodology relies on the assumption of the conditional mean-scale model (\ref{1.1}). This is a limitation compared to more general model structures that only assume $E[\varepsilon|X] = 0$. It is interesting to develop goodness-of-fit tests for more general regression models when the dimension of predictor vector is divergent or even larger than the sample size. The research is ongoing.

\section*{Acknowledgement}
We thank the Editor, the Associate Editor and two referees for their constructive comments and suggestions, which greatly improved the quality of this paper. The thank also goes to Miss Jiaqi Huang for her effort in the presentation of this paper.
Falong Tan was supported by the National Natural Science Foundation of China (12071119). Xu Guo was supported by the National Natural Science Foundation of China (12071038). Lixing Zhu was supported by National Natural Science Foundation of China (12131006, 12471276)).

\section*{Supplementary material}
Supplementary Material available online includes all proofs for the theoretical results presented in the main text, and the tests for checking the parametric conditional variance function based on the residual empirical process and its martingale transformation in diverging-dimension settings. It also contains typical examples for illustrating the technical conditions and theoretical results in the main text.

\bibliographystyle{apalike}
\bibliography{martingale}

\end{document}